\documentclass[trsc,nonblindrev]{informs3} 

\OneAndAHalfSpacedXI 

\usepackage{natbib}
 \bibpunct[, ]{(}{)}{,}{a}{}{,}%
 \def\bibfont{\small}%
\TheoremsNumberedThrough     

\EquationsNumberedThrough    

\usepackage{mathtools}
\usepackage{multirow}
\usepackage{tikz}
\usepackage{subcaption}
\usepackage{comment}
\usepackage{accents} 
\usepackage{algorithmic} 

\usepackage[labelformat=simple]{subcaption}

\usepackage{bm}
\newcommand{\vc}[1]{\bm{#1}}
\newcommand{\mr}[1]{\mathrm{#1}}
\newcommand{\prob}{{\text{Pr}}}

\def\mc{\mathcal}
\def\Z{\mathbb{Z}}
\def\Re{{\mathbb R}}

\def\Sol{\mr{Sol}}

\newcommand{\blue}[1]{\textcolor{blue}{#1}}

\newcommand{\suchthat}{\;\ifnum\currentgrouptype=16 \middle\fi|\;}
\newcommand{\sign}{\operatorname{sign}}
\newcommand{\ind}{\operatorname{ind}}

\usepackage[ruled,vlined,linesnumbered]{algorithm2e}

\begin{document}


\RUNAUTHOR{Salemi and Davarnia}

\RUNTITLE{Solving Unsplittable Network Flow Problems with Decision Diagrams}

\TITLE{Solving Unsplittable Network Flow Problems with Decision Diagrams}

\ARTICLEAUTHORS{%
\AUTHOR{Hosseinali Salemi, Danial Davarnia}
\AFF{Department of Industrial and Manufacturing Systems Engineering, Iowa State University, Ames, IA 50011, \EMAIL{hsalemi@iastate.edu}, \EMAIL{davarnia@iastate.edu}}
} 

\ABSTRACT{%
In unsplittable network flow problems, certain nodes must satisfy a combinatorial requirement that the incoming arc flows cannot be split or merged when routed through outgoing arcs. 
This so-called \emph{no-split no-merge} requirement arises in unit train scheduling where train consists should remain intact at stations that lack necessary equipment and manpower to attach/detach them.
Solving the unsplittable network flow problems with standard mixed-integer programming formulations is computationally difficult due to the large number of binary variables needed to determine matching pairs between incoming and outgoing arcs of nodes with no-split no-merge constraint.
In this paper, we study a stochastic variant of the unit train scheduling problem where the demand is uncertain.
We develop a novel decision diagram (DD)-based framework that decomposes the underlying two-stage formulation into a master problem that contains the combinatorial requirements, and a subproblem that models a continuous network flow problem.
The master problem is modeled by a DD in a transformed space of variables with a smaller dimension, leading to a substantial improvement in solution time.
Similarly to the Benders decomposition technique, the subproblems output cutting planes that are used to refine the master DD.
Computational experiments show a significant improvement in solution time of the DD framework compared with that of standard methods.              
}%


\HISTORY{}

\maketitle

%


\section{Introduction} \label{sec: Introduction}
\begin{figure}[!hbt]
\centering
\includegraphics[scale=0.5]{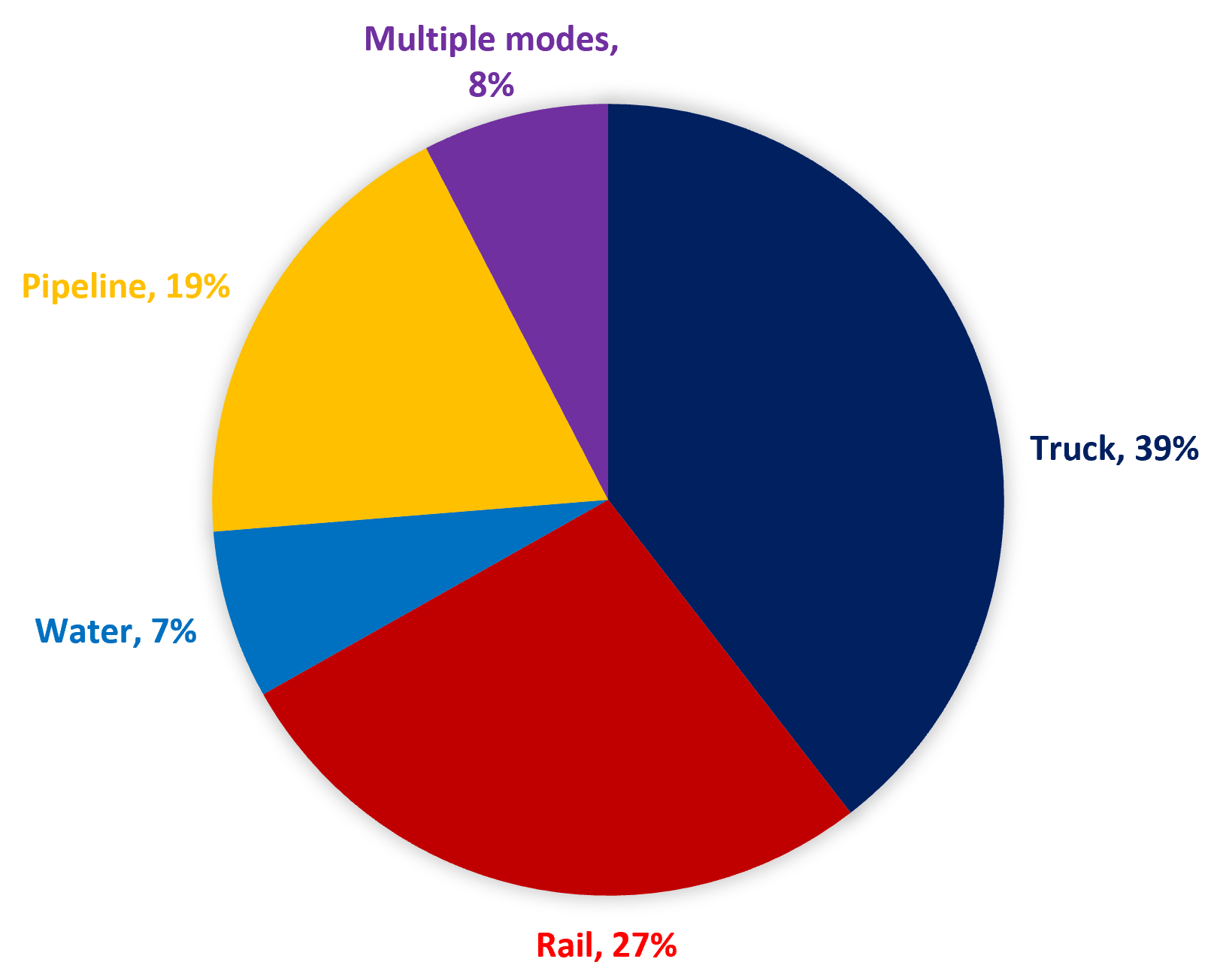}
\caption{Pie chart for ton-miles of freight shipments by mode within the U.S.\ in 2018. Multiple modes includes mail. Air and truck-air with the share of $0.1\%$ are omitted.}
\label{fig:freight_shipments}
\end{figure}

Over the past several decades, rail freight transportation has continued to grow as the prime means of transportation for high-volume commodities.
Advantages of rail transportation include reliability, safety, cost-efficiency and environmental-sustainability as compared with alternative methods of transportation.
In terms of scale, the rail network accounted for $27.2$ percent of U.S.\ freight shipment by ton-miles in 2018~\citep{furchtgott2021pocket}; see Figure~\ref{fig:freight_shipments}. 
The Federal Highway Administration estimates that the total U.S.\ freight shipments will be 24.1 billion tons in 2040, a $30$ percent increase from the 2018 total transportation of 18.6 billion tons. 
With the purpose of meeting such market growth, America's freight railway companies have invested nearly $\$740$ billion on capital expenditures and maintenance from 1980 to 2020~\citep{AAR}.      

To reduce rail freight transportation costs and shipment delays, railroad companies offer \emph{unit train} services for carrying high-volume products. Unit trains haul a single type freight in a way that no car is attached or detached while the cargo train is on its way from an origin to a destination, except in specific locations that are equipped with required manpower and machinery. 
These trains usually operate all day, use dedicated equipment, and can be loaded/unloaded in 24 hours. They are known to be one of the fastest and most efficient means of railroad transportation.~\citep{AAR}.
Traditionally, unit trains are used to carry bulk cargo such as coal, grain, cement, and rock. 
Bulk liquids like crude oil and food such as wheat and corn are also shipped by unit trains.
According to the Federal Railroad Administration data, bulk commodities account for $91$ percent of the U.S.\ railroad freights. 
Approximately all coal shipped through railways in the U.S. are transported by unit trains. Moreover, these trains contribute significantly to the shipping process of crude oil as each unit train is capable of carrying 85,000 barrels~\citep{AAR}.
In an operational level, the core unit train model can be described as follows. Given a set of supply, intermediate, and demand locations in a railroad network, the unit train scheduling problem seeks to find optimal routes for unit trains to send flows from supply to demand points with the objective of minimizing the total transportation cost while meeting demand of customers, respecting capacities of tracks, and satisfying no-car attaching/detaching requirements in specific locations.
As a result, designing \emph{blocking plans} to determine locations where cars need to be switched between trains is irrelevant in this problem, unlike scheduling other types of trains~\citep{davarnia2019network}.

Despite the significance of unit train scheduling, exact optimization approaches to solve associated problems are scarce, partially due to their structural complexities. 
One of the main challenges in modeling unit trains is the requirement that the train consists must remain intact when passing through stations that lack necessary busting/formation equipment.
In optimization, this requirement is referred to as \textit{no-split no-merge} (NSNM), which guarantees that the flows entering to or exiting from certain nodes of the unit train network cannot be split or merged.
Incorporating this requirement into typical transportation network models yields the so-called \textit{generalized unsplittable flow problem} (GUFP), where the objective is to determine the minimum-cost unit train schedules that satisfy the given demand. 
Numerous studies have shown that considering deterministic demands might result in the complete failure of the transportation scheduling~\citep{demir2016green, layeb2018simulation}, motivating the study of stochastic variants of the unit train scheduling problems where the demand is uncertain. 
As a result, in this paper, we consider a stochastic variant of the GUFP, referred to SGUFP, that is modeled as a two-stage optimization problem.
The first stage decides a matching between the incoming and outgoing arcs of the nodes of the railroad network, and the second stage determines the amount of flow that should be sent through the matching arcs of the network to satisfy the uncertain demand represented by a number of demand scenarios. We propose a novel exact solution framework to solve this problem in the operational level.

Our proposed methodology is based on \textit{decision diagrams} (DDs), which are compact graphical data structures. 
DDs were initially introduced to represent boolean functions with applications in circuit design. 
Over the past decade, researchers have successfully extended DDs domain by developing DD-based algorithms to solve discrete optimization problems in different areas of application.
Because of its structural limitation to model integer programs only, DDs have never been used to solve transportation problems that inherently include continuous variables.
In this paper, we extend the application scope of DDs by introducing a novel framework that is capable of modeling network problems with both integer and continuous components as in the SGUFP.

\subsection{Literature Review on Train Scheduling} \label{subsec: Background on Train Scheduling}
Many variants of train routing and scheduling problems with different objective functions and set of constraints under deterministic and stochastic conditions have been introduced and vastly studied in the literature; see surveys by~\cite{cordeau1998survey},~\cite{harrod2010operations},~\cite{lusby2011railway},~\cite{cacchiani2012nominal}, and~\cite{turner2016review} for different problems classifications and structures. 
Mixed integer linear and nonlinear programming formulations are among the most frequent exact approaches to model different classes of these problems~\citep{jovanovic1991tactical, huntley1995freight,sherali1998tactical,lawley2008time,haahr2017integrating,davarnia2019network}. 
Proposed solution techniques include but are not limited to branch-and-bound methods~\citep{jovanovic1991tactical, fuchsberger2007solving}, branch-and-cut frameworks~\citep{zwaneveld2001routing,ceselli2008optimizing}, branch-and-price approaches~\citep{lusby2008optimization,lin2016branch}, graph coloring algorithms~\citep{cornelsen2007track}, and heuristics~\citep{carey2007scheduling,liu2011optimising,iccyuz2016two}.     
Rolling stock scheduling~\citep{abbink2004allocation, alfieri2006efficient, haahr2016comparison, borndorfer2016integrated} that assigns rolling stocks to a given timetable, and crew scheduling~\citep{kwan2011case, shen2013evolutionary, heil2020railway} that covers train activities by assigning crews to the associated operations are other major problems arising in the area of railroad planning.

Due to the inherent uncertainty in different types of train scheduling and routing problems, many researchers have studied stochastic variants of the problems where the supply/demand is considered to be uncertain. \cite{jordan1983stochastic} propose a model for railroad car distribution where supply and demand of cars are uncertain.
\cite{jin2019approach} study a chance-constrained programming model for the train stop planning problem under stochastic demand. \cite{ying2020actor} propose a deep reinforcement learning approach for train scheduling where the passenger demand is uncertain. Recently, \cite{gong2021train} propose a stochastic optimization method to solve a train timetabling problem with uncertain passenger demand. Also see works by~\cite{meng2011robust, quaglietta2013stability,larsen2014susceptibility} that consider train dispatching problems under stochastic environments.

In the context of unit train scheduling, \cite{lawley2008time} study a time-space network flow model to schedule bulk railroad deliveries for unit trains. 
In their model, the authors consider characteristics of underlying rail network, demands of customers, and capacities of tracks, stations, and loading/unloading requirements. 
They propose a mixed integer programming (MIP) formulation that maximizes the demand satisfaction while minimizing the waiting time at stations. \cite{lin2014two} (cf.~\cite{lin2016branch}) propose a model for a train scheduling problem that is capable to capture locations where coupling/decoupling is forbidden. 
They develop a branch-and-price algorithm inspired by column generation to solve the associated problem. \cite{lin2018redundant} also propose a heuristic branch-and-bound approach to decrease coupling/decoupling redundancy. \cite{iccyuz2016two} study the problem of planning coal unit trains that includes train formation, routing, and scheduling. 
As noted by the authors, their proposed MIP formulation fails to solve the problem directly due to its large size.
As a remedy, they develop a time-efficient heuristic that produces good quality solutions. 
More recently, \cite{davarnia2019network} introduce and study the GUFP with application to unit train scheduling. 
In particular, the authors show how to impose NSNM restrictions in network optimization problems. 
They present a polyhedral study and propose a MIP formulation to model a stylized variant of the unit train scheduling problem. 
In the present paper, we use their formulation (see section~\ref{subsec: MIP Formulation}) as a basis for our solution framework.    

The unsplittable flow problem (UFP) was first introduced by~\cite{kleinberg1996approximation} as a generalization of the disjoint path problem. 
Given a network with capacities for arcs and a set of source-terminal vertex pairs with associated demands and rewards, the objective in the UFP is to maximize the total revenue by selecting a subset of source-terminal pairs and routing flows through a \emph{single} path for each of them to satisfy the associated demand. 
In the GUFP, however, there can exist nodes that do not need to respect the NSNM requirement, and demands can be satisfied by passing flows through multiple paths. 
It is well-known that different variants of UFP are NP-hard~\citep{baier2005k,kolman2006improved,chakrabarti2007approximation}. 
Since its introduction, the UFP structure has been used in different areas of application, from bandwidth allocation in heterogeneous networks~\citep{kolman2006improved}, to survivable connection-oriented networks~\citep{walkowiak2006new}, and virtual circuit
routing problems~\citep{hu2009algorithm}. 
Considering the hardness of the problem, approximation algorithms have been a common technique to tackle different variants of the UFP in the literature~\citep{baier2005k,chakrabarti2007approximation}. 

\subsection{Literature Review on Decision Diagrams} \label{subsec: Background on Decision Diagrams}
DDs are directed acyclic graphs with a source and a terminal node where each source-terminal path encodes a feasible solution to an optimization problem. 
In DDs, each layer from the source to the terminal represents a decision variable where labels of arcs show their values.
\cite{hadzic2006discrete} proposed to use DDs to model the feasible region of a discrete optimization problem and used it for postoptimality analysis. 
Later,~\cite{andersen2007constraint} presented relaxed DDs to circumvent the exponential growth rate in the DD size when modeling large discrete optimization problems. 
\cite{bergman2016discrete} introduced a branch-and-bound algorithm that iteratively uses relaxed and restricted DDs to find optimal solution. 
The literature contains many successful utilization of DDs in different domains; see works by~\cite{bergman2018discrete, serra2019compact, davarnia2020outer, gonzalez2020integrated}, and~\cite{hosseininasab2021exact} for some examples. 

Until recently, applications of DDs were limited to discrete problems, and the question on how to use DDs in solving optimization problems with continuous variables was unanswered. 
To address this limitation, \cite{davarnia2021strong} proposed a technique called arc-reduction that generates a DD that represents a relaxation of the underlying continuous problem. 
In a follow-up work, \cite{salemistructure} established necessary and sufficient conditions for a general MIP to be representable by DDs. 
They showed that a bounded MIP can be remodeled and solved with DDs through employing a specialized Benders decomposition technique. 
In this paper, we build on this framework to design a novel DD-based methodology to solve the SGUFP.             

\subsection{Contributions} \label{subsec: our contributions}
While there are several studies in the literature dedicated to the unit train problem, exact methodologies that provide a rigorous treatment of the NSNM requirement at the heart of unit train models are scarce.
In this paper, we design a novel exact DD-based framework to solve the SGUFP, as a more realistic and more challenging variant of this problem class. 
To our knowledge, this is the first work that studies SGUFP from an exact perspective, and the first application of DDs to a transportation problem.
Our proposed framework formulates the problem in a transformed space of variables, which has a smaller dimension compared to the standard MIP formulations of the SGUFP.
This presentation mitigates the computational difficulties stemmed from the MIP formulation size, providing a viable solution approach for large-scale network problems.
The core principles of our DD framework can also be used to model other transportation problems with similar structure, as an alternative to traditional network optimization techniques. 

The remainder of this paper is organized as follows. 
In Section~\ref{sec: Background on DDs} we provide basic definitions and a brief overview on discrete and continuous DD models, including the DD-BD method to solve bounded MIPs. 
In Section~\ref{sec: DDBD}, we adapt the DD-BD method to solve the SGUFP. 
We propose algorithms to construct exact and relaxed DDs to solve the problem in a transformed space. 
Section~\ref{sec: Computational Experiments} presents computational experiments to evaluate the performance of the DD-BD method for the SGUFP. 
We give concluding remarks in Section~\ref{sec: Conclusion}.

\section{Background on DDs} \label{sec: Background on DDs}
In this section, we present basic definitions and results relevant to our DD analysis.

\subsection{Overview}
A DD $\mc D = (\mc U,\mc A, l)$ with node set $\mc U$, arc set $\mc A$, and arc label mapping $l: \mc A \to \Re$ is a directed acyclic graph with $n \in \mathbb{N}$ arc layers $\mc A_1, \mc A_2, \dots, \mc A_n$, and $n+1$ node layers $\mc U_1,\mc U_2,\dots,\mc U_{n+1}$.
The node layers $\mc U_1$ and $\mc U_{n+1}$, with $|\mc U_1| = |\mc U_{n+1}| = 1$, contain the root $r$ and the terminal $t$, respectively. 
In any arc layer $j \in [n] \coloneqq \{1,2,\dots, n\}$, an arc $(u,v) \in \mc A_j$ is directed from the tail node $u \in \mc U_j$ to the head node $v \in \mc U_{j+1}$. 
The \textit{width} of $\mc D$ is defined as the size of its largest $\mc U_j$. 
DDs can model a bounded integer set $\mc P \subseteq \Z^n$ in such a way that each $r$-$t$ arc-sequence (path) of the form $(a_1, \dotsc, a_n) \in \mc A_1 \times \dotsc \times \mc A_n$ encodes a point $\vc y \in \mc P$ where $l(a_j) = y_j$ for $j \in [n]$, that is $\vc y$ is an $n$-dimensional point in $\mc P$ whose $j$-th coordinate is equal to the label value $l(a_j)$ of arc the $a_j$.
For such a DD, we have $\mc P = \Sol(\mc D)$, where $\Sol(\mc D)$ denotes the finite collection of all $r$-$t$ paths.         

The graphical property of DDs can be exploited to optimize an objective function over a discrete set $\mc{P}$. 
To this end, DD arcs are weighted in such a way that the cumulative weight of an $r$-$t$ path that encodes a solution $\vc y \in \mc P$ equals to the objective function value evaluated at $\vc y$. 
Then, a shortest (resp. longest) $r$-$t$ path for the underlying minimization (resp. maximization) problem is found, an operation that can be performed in polynomial time. 

The construction of an \textit{exact} DD as described above is computationally prohibitive due to the exponential growth rate of its size.
To alleviate this difficulty, \emph{relaxed} and \emph{restricted} DDs are proposed to keep the size of DDs under control. 
In a relaxed DD, nodes are merged in such a way that the width of the resulting diagram is bounded by a predetermined width limit. 
This node-merging process ensures that all feasible solutions of the original set are encoded by a subset of all $r$-$t$ paths in the resulting DD.
Optimization over this relaxed DD provides a dual bound to the optimal solution of the original problems.
In a restricted DD, the collection of all $r$-$t$ paths of the DD encode a subset of the feasible solutions of the original set.
Optimization over this restricted DD provides a primal bound to the optimal solution of the original problems.
The restricted and relaxed DDs can be iteratively refined in a branch-and-bound scheme to find the optimal value of a problem through convergence of their primal and dual bounds.
The following example illustrates an exact, relaxed and restricted DD for a discrete optimization problem.     
\begin{example}\label{ex: IP} 
Consider the discrete optimization problem $\max\{5y_1 + 10y_2 + 4y_3 \suchthat \vc y \in \mc P\}$ where $\mc P = \{(1,0,0),(1,0,1),(0,1,0),(0,0,1),(0,0,0)\}$.
The exact DD $\mc D$ with width 3 in Figure~\ref{subfig: exact DD} models the feasible region $\mc P$. The weight of each arc $a \in \mc A_j$, for $j \in \{1,2,3\}$, shows the contribution of variable $y_j$'s value assignment to the objective function. 
The longest $r$-$t$ path that encodes the optimal solution $(y^*_1,y^*_2,y^*_3)=(0,1,0)$ has length 10, which is the optimal value to the problem. 
By reducing the width limit to 2, we can build relaxed and restricted DDs for $\mc{P}$ as follows.
The relaxed DD $\overline{\mc D}$ in Figure~\ref{subfig: relaxed DD} provides an upper bound to the optimal solution, where the longest path with length 14 is obtained by an infeasible point $(\overline{y}_1,\overline{y}_2,\overline{y}_3)=(0,1,1)$. 
Finally, the restricted DD $\underline{\mc D}$ in Figure~\ref{subfig: restricted DD} gives a lower bound to the optimal solution, where the longest path with length 9 encodes a feasible solution $(\underline{y_1},\underline{y_2},\underline{y_3})=(1,0,1)$.     

\begin{figure}[!htbp]
\captionsetup[subfigure]{justification=centering}
\begin{subfigure}[b]{0.3\linewidth}
\centering
\begin{tikzpicture}[scale=0.8]
\begin{scope}[every node/.style={circle,color=blue,thick,draw,minimum size=0.4cm}]
\node(1) at (0,0) {$r$};
\node(2) at (1.4,-2) {};
\node(3) at (-1.4,-2) {};
\node(4) at (2.9,-4) {};
\node(5) at (0,-4) {};
\node(6) at (-2.9,-4) {};
\node(7) at (0,-6.5) {$t$};
\end{scope}
\path (1) edge node [right] {\small 5} (2);
\path[color=purple,dashed] (1) edge node [left] {\small 0} (3);
\path[dashed] (2) edge node [right] {\small 0} (4);
\path[color=purple] (3) edge node [right] {\small 10} (5);
\path[dashed] (3) edge node [left] {\small 0} (6);
\path[dashed] (4) edge[bend left=20] node [right] {\small 0} (7);
\path(4) edge[bend right=20] node [above] {\small 4} (7);
\path[color=purple,dashed] (5) edge node [right] {\small 0} (7);
\path(6) edge[bend left=20] node [above] {\small 4} (7);
\path[dashed] (6) edge[bend right=20] node [left] {\small 0} (7);
%
%
\end{tikzpicture}
\caption{\small Exact DD $\mc D$}
\label{subfig: exact DD}
\end{subfigure}\hspace{0.05\linewidth}
\begin{subfigure}[b]{0.3\linewidth}
\centering
\hspace{0.2cm}
\begin{tikzpicture}[scale=0.8]
\begin{scope}[every node/.style={circle,color=blue,thick,draw,minimum size=0.4cm}]
\node(1) at (0,0) {$r$};
\node(2) at (1.4,-2) {};
\node(3) at (-1.4,-2) {};
\node(4) at (1.4,-4) {};
\node(5) at (-1.4,-4) {};
\node(6) at (0,-6.5) {$t$};
\end{scope}
\path (1) edge node [right] {\small 5} (2);
\path[color=purple,dashed] (1) edge node [left] {\small 0} (3);
\path[dashed] (2) edge node [right] {\small 0} (4);
\path[color=purple] (3) edge node [above] {\small 10} (4);
\path[dashed] (3) edge node [left] {\small 0} (5);
\path[dashed] (4) edge[bend left=20] node [right] {\small 0} (6);
\path[color=purple] (4) edge[bend right=20] node [left] {\small 4} (6);
\path (5) edge[bend left=20] node [right] {\small 4} (6);
\path[dashed] (5) edge[bend right=20] node [left] {\small 0} (6);
%
%
\end{tikzpicture}
\caption{\small Relaxed DD $\overline{\mc D}$}
\label{subfig: relaxed DD}
\end{subfigure}\hspace{0.05\linewidth}
\begin{subfigure}[b]{0.3\linewidth}
\centering
\hspace{0.6cm}
\begin{tikzpicture}[scale=0.8]
\begin{scope}[every node/.style={circle,color=blue,thick,draw,minimum size=0.4cm}]
\node(1) at (0,0) {$r$};
\node(2) at (1.4,-2) {};
\node(3) at (-1.4,-2) {};
\node(4) at (1.4,-4) {};
\node(5) at (-1.4,-4) {};
\node(6) at (0,-6.5) {$t$};
\end{scope}
\path[color=purple] (1) edge node [right] {\small 5} (2);
\path[dashed] (1) edge node [left] {\small 0} (3);
\path[color=purple][dashed] (2) edge node [right] {\small 0} (4);
\path[dashed](3) edge node [left] {\small 0} (5);
\path[dashed] (4) edge[bend left=20] node [right] {\small 0} (6);
\path[color=purple] (4) edge[bend right=20] node [left] {\small 4} (6);
\path (5) edge[bend left=20] node [right] {\small 4} (6);
\path[dashed] (5) edge[bend right=20] node [left] {\small 0} (6);
\hspace{0.5cm}\begin{scope}[every node/.style={circle,color=black,minimum size=0.1cm}]
\node(6) at (2.2,-1)  {$y_1$};
\node(7) at (2.2,-3.2)  {$y_2$};
\node(8) at (2.2,-5.2)  {$y_3$};
\end{scope}
\end{tikzpicture}
\caption{\small Restricted DD $\underline{\mc D}$}
\label{subfig: restricted DD}
\end{subfigure}
\caption{The exact, relaxed, and restricted DDs representing $\mc P$ in Example~\ref{ex: IP}. Solid and dotted arcs indicate one and zero arc labels, respectively. Numbers next to arcs represent weights.}
\end{figure}
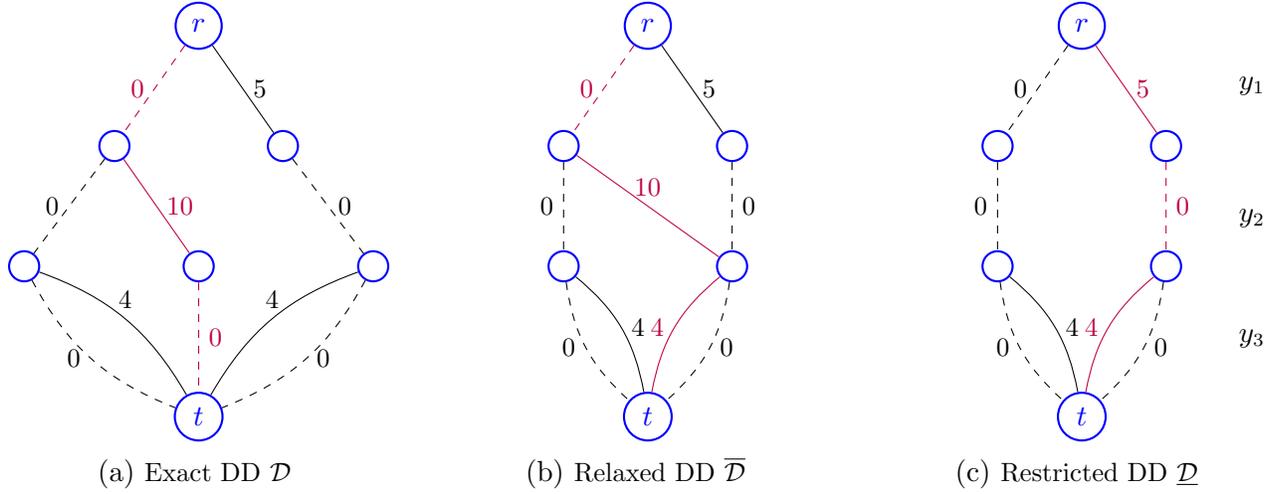

\end{example}

\subsection{Continuous DD Models} \label{subsec: Continuous DD Models}
While the framework described in the previous section can be applied to solve different classes of discrete optimization problems, its extension to model sets with continuous variables requires a fundamentally different approach. 
The reason that the traditional DD structure is not viable for continuous sets is that representing the domain of a continuous variable through arcs requires an infinite number of them, spanning all values within a continuous interval, which is structurally prohibitive in DD graphs. 
Fortunately, there is a way to overcome this obstacle by decomposing the underlying set into certain rectangular formations, which can in turn be represented through node-sequences in DDs.
In what follows, we give an overview of these results as relevant to our analysis.

Consider a bounded set $\mc P \subseteq \Re^n$.
\cite{salemistructure} give necessary and sufficient conditions for $\mc P$ to admit the desired rectangular decomposition.
Such a set is said to be \textit{DD-representable} w.r.t.\ a fixed index set $I \subseteq [n]$, as there exists a DD $\mc D$ such that $\max\{f(\vc x) \suchthat \vc x \in \mc P\} = \max\{f(\vc x) \suchthat \vc x \in \Sol(\mc D)\}$ for every function $f(\vc x)$ that is convex in the space of variables $\vc x_I$. 
A special case of DD-representable sets is given next.


\begin{proposition} \label{prop: mip}
Any bounded mixed integer set of the form $\mc P \subseteq \Z^n \times \Re$ is DD-representable w.r.t. $I=\{n+1\}$. 
\Halmos
\end{proposition}

This result gives rise to a novel DD-based framework to solve general bounded MIPs as outlined below.
Consider a bounded MIP $\mc H \coloneqq \max\{\vc c \vc y + \vc d \vc x \suchthat A \vc y + G \vc x \le \vc b,~\vc y\in \Z^n\}$. 
Using Benders decomposition (BD), formulation $\mc H$ is equivalent to $\max_{\vc y \in \Z^n}\{\vc{cy} + \max_{\vc x}\{\vc{dx} \suchthat G \vc x \le \vc b - A \vc y\}\}$, which can be reformulated as $\mc M = \max\{\vc{cy} + z \suchthat (\vc y;z)\in \Z^n \times [l,u]\}$, where $l, u \in \Re$ are some valid bounds on $z$ induced from the boundedness of $\mc H$. 
Here, $\mc M$ is the master problem and $z$ represents the objective value of the subproblem $\max_{\vc x}\{\vc{dx} \suchthat G \vc x \le \vc b - A \bar{\vc{y}}\}$ for any given $\bar{\vc y}$ as an optimal solution of the master problem.
The outcome of the subproblems is either an optimality cut or a feasibility cut that will be added to the master problem. Then, the master problem will be resolved.
Proposition~\ref{prop: mip} implies that formulation $\mc M$ can be directly modeled and solved with DDs. 
For this DD, we assign $n$ arc layers to the integer variables $y_1,y_2,\dots,y_n$, and one arc layer to the continuous variable $z$ with only two arc labels showing a lower and upper bound for this variable.
To find an optimal solution, the longest path is calculated, which will be used to solve the subproblems.
Note that since $\mc M$ is a maximization problem, a longest path of the associated DD encodes an optimal solution, and its length gives the optimal value; see Example~\ref{ex: MIP}.
The feasibility and optimality cuts generated by the subproblems will then be added to \textit{refine} the DD, whose longest path will be recalculated.
The refinement technique consists of removing arcs of the DD that lead to solutions that violate the added inequality, as well as splitting nodes of the DD that lead to different subsequent partial assignments; see \cite{bergman2016decision} for a detailed account on DD refinement techniques.
We illustrate this approach in Example~\ref{ex: MIP}.

\begin{example} \label{ex: MIP} 
Suppose that $\max\{2y_1 + 4y_2 + z \suchthat \vc y \in \mc P, z \le 25\}$ forms the master problem at the penultimate iteration of a BD algorithm, where $\mc P =\{(0,0),(1,1)\}$.
This problem is represented by the DD $\mc D$ in Figure~\ref{subfig: penultimate iteration} where $-M$ is a valid lower bound for $z$.
The longest path of $\mc D$ encodes the solution $(\hat y_1, \hat y_2, \hat z) = (1,1,25)$. 
Assume that using the point $(\hat y_1,\hat y_2)=(1,1)$ in the associated subproblem generates an optimality cut $z \le 3y_1 + 2y_2 + 10$ for the final iteration of the BD algorithm. 
Refining DD $\mc D$ with respect to this cut yields the new DD in Figure~\ref{subfig: final iteration}. 
The longest path represents the optimal solution $(y^*_1,y^*_2,z^*)=(1,1,15)$ with length 21, which is the optimal value.     

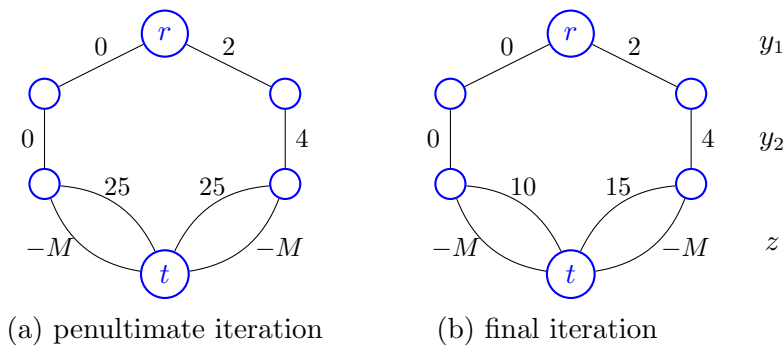
\begin{figure}[!htbp]
\centering
\captionsetup[subfigure]{justification=centering}
\begin{subfigure}[b]{0.3\linewidth}
\centering
\begin{tikzpicture}[scale=0.8]
\begin{scope}[every node/.style={circle,color=blue,thick,draw,minimum size=0.4cm}]
\node(1) at (0,0) {$r$};
\node(2) at (2,-1) {};
\node(3) at (-2,-1) {};
\node(4) at (2,-2.5) {};
\node(5) at (-2,-2.5) {};
\node(6) at (0,-4) {$t$};
\end{scope}
\path(1) edge node [above] {\small 2} (2);
\path(1) edge node [above] {\small 0} (3);
\path(2) edge node [right] {\small 4} (4);
\path(3) edge node [left] {\small 0} (5);
\path(4) edge[bend left=30] node [right] {\small $-M$} (6);
\path(4) edge[bend right=30] node [above] {\small 25} (6);
\path(5) edge[bend left=30] node [above] {\small 25} (6);
\path(5) edge[bend right=30] node [left] {\small $-M$} (6);
\end{tikzpicture}
\caption{penultimate iteration}
\label{subfig: penultimate iteration}
\end{subfigure}
\begin{subfigure}[b]{0.3\linewidth}
\centering
\hspace{0.5cm}
\begin{tikzpicture}[scale=0.8]
\begin{scope}[every node/.style={circle,color=blue,thick,draw,minimum size=0.4cm}]
\node(1) at (0,0) {$r$};
\node(2) at (2,-1) {};
\node(3) at (-2,-1) {};
\node(4) at (2,-2.5) {};
\node(5) at (-2,-2.5) {};
\node(6) at (0,-4) {$t$};
\end{scope}
\path(1) edge node [above] {\small 2} (2);
\path(1) edge node [above] {\small 0} (3);
\path(2) edge node [right] {\small 4} (4);
\path(3) edge node [left] {\small 0} (5);
\path(4) edge[bend left=30] node [right] {\small $-M$} (6);
\path(4) edge[bend right=30] node [above] {\small 15} (6);
\path(5) edge[bend left=30] node [above] {\small 10} (6);
\path(5) edge[bend right=30] node [left] {\small $-M$} (6);
\hspace{1.4cm}\begin{scope}[every node/.style={circle,color=black,minimum size=0.1cm}]
\node(6) at (1.6,-0.2)  {$y_1$};
\node(7) at (1.6,-1.8)  {$y_2$};
\node(8) at (1.6,-3.5)  {$z$};
\end{scope}
\end{tikzpicture}
\caption{final iteration}
\label{subfig: final iteration}
\end{subfigure}
\caption{The last two iterations of solving the master problem in Example~\ref{ex: MIP}}
\end{figure}

\end{example}

Using the DD framework as outlined above can be computationally challenging due to exponential growth rate of the size of an exact DD. 
To mitigate this difficulty, restricted/relaxed DDs can be employed inside of the BD framework as demonstrated in Algorithm~\ref{alg: DD-BD}. 
We refer to this solution method as DD-BD \citep{salemistructure}.

\begin{algorithm}[!ht]
\caption{DD-BD}
\label{alg: DD-BD}
\KwData{MIP $\mc H$, construction method to build restricted and relaxed DDs for $\mc M$}
\KwResult{An optimal solution $(\vc y^*,z^*)$ and optimal value $w^*$ to $\mc H$}

initialize set of partial assignments $\mc {\hat Y} = \{\ominus\}$, set of Benders cuts $\mc C=\emptyset$, and $w^* = -\infty$ 

\eIf{$\mc {\hat Y} =  \emptyset$}
{terminate and return $(\vc{y^*},z^*)$ and $w^*$}
{
select $\hat{\vc y} \in \hat{\mc Y}$ and update $\mc{\hat Y} \leftarrow \mc{\hat Y} \setminus \{\hat{\vc y}\}$

create a restricted DD $\underline{\mc D}$ associated with $\mc M^C(\hat{\vc y})$

\eIf{$\underline{\mc D} \neq \emptyset$}
{
find a longest $r$-$t$ path of $\underline{\mc D}$ with encoding point $(\underline{\vc y}, \underline{z})$ and length $\underline{w}$ 

solve the BD subproblem using $\underline{\vc y}$ to obtain Benders cut $\underline{C}$

\eIf{$\underline{C} \in \mc C$}
{go to line 17}
{
update $\mc C \leftarrow \mc C \cup \underline{C}$ and refine $\underline{\mc D}$ w.r.t. $\underline{C}$

go to line 8
}
}
{go to line 2}

\If{$\underline{w} > w^*$}{update $w^* \leftarrow \underline{w}$ and $(\vc y^*,z^*) \leftarrow (\underline{\vc y}, \underline{z})$}

\eIf{$\underline{\mc D}$ provides an exact representation of $\mc M^C(\hat{\vc y})$}{go to line 2}
{
create a relaxed DD $\overline{\mc D}$ associated with $\mc M^C(\hat{\vc y})$

find a longest $r$-$t$ path of $\overline {\mc D}$ with length $\overline{w}$

\If{$\overline{w} > w^*$}
{

solve the BD subproblem using $\overline{\vc y}$ to obtain Benders cut $\overline{C}$

\eIf{$\overline{C} \in \mc C$}
{go to line 31}
{
update $\mc C \leftarrow \mc C \cup \overline{C}$ and refine $\overline{\mc D}$ w.r.t. $\overline{C}$

go to line 23
}

\ForAll{$u$ in the last exact layer of $\overline{\mc D}$}
{update $\hat{\mc Y} \leftarrow \hat{\mc Y} \cup \{\tilde{\vc y}\}$ where $\tilde{\vc y}$ encodes longest $r$-$u$ path of $\overline{\mc D}$}
}
}

go to line 2
}

\end{algorithm}

In explaining the steps of Algorithm~\ref{alg: DD-BD}, let point $\hat{\vc y} \in \Z^k$, where $k \le n$, be a partial value assignment to the first $k$ coordinates of variable $\vc y$, i.e., $y_i = \hat y_i$ for all $i \in [k]$. 
We record the set of all partial value assignments in $\hat {\mc{Y}} = \{\hat{\vc y} \in \Z^k \suchthat k \in [n]\} \cup \{\ominus\}$, where $\ominus$ represents the case where no coordinate of $\vc y$ is fixed. 
Set $\mc C$ contains the produced Benders cuts throughout the algorithm, and we denote the feasible region described by these cuts by $\mc F^{\mc C}$. 
Further, define $\mc M^C(\hat{\vc y})= \max\{\vc{cy} + z \suchthat (\vc y; z) \in \Z^n \times [l,u] \cap \mc F^{\mc C},~y_i = \hat y_i,\forall i \in [k]\}$ to be the restricted master problem $\mc M$ obtained through adding cuts in $\mc C$ and fixing the partial assignment $\hat{\vc y}$. 
In this definition, the case with $\mc C = \emptyset$ and $\hat{\mc Y} = \{\ominus\}$ is denoted by $\mc M^{\emptyset}(\ominus) = \mc M$, which is an input to Algorithm~\ref{alg: DD-BD}.  

The algorithm starts with constructing a restricted DD $\underline{\mc D}$ corresponding to $\mc M^C(\hat{\vc y})$ with empty initial values for $C$ and $\hat{y}$.
We then find a longest $r$-$t$ path of $\underline{\mc D}$ encoding solution $(\underline{\vc y},\underline{z})$. 
Next, using $\underline{\vc y}$, we solve the associated subproblem to obtain a feasibility/optimality cut $\underline{C}$.
We add this cut to $\mc C$, refine $\underline{\mc D}$ according to it, and find a new longest $r$-$t$ path. 
We repeat these steps until no new feasibility/optimality cut is generated. 
At this point, the length of a longest $r$-$t$ path of $\underline{\mc D}$, denoted by $\underline{w}$, gives a lower bound to the master problem $\mc M$, which is also a valid lower bound to the original problem $\mc H$. 
The value of $\underline{w}$ can be used to update $w^*$, the optimal value of $\mc H$ at termination. Next, we create a relaxed DD $\overline{\mc D}$ corresponding to $\mc M^C(\hat{\vc y})$. 
We find a longest $r$-$t$ path of $\overline{\mc D}$ that provides an upper bound $\overline{w}$ to $\mc M$. 
If the upper bound $\overline{w}$ is strictly greater than the current value of $w^*$, we follow steps similarly to the case for $\underline{\mc D}$ to iteratively refine $\overline{\mc D}$ w.r.t. feasibility/optimality cuts through solving the subproblems, until no new cut is generated.
Next, we perform a specialized branch-and-bound procedure to improve the bound through expanding merged layers of the DD. 
To this end, we add all the partial assignments associated with nodes in the last exact layer of $\overline{D}$ (the last node layer in which no nodes are merged) to the collection $\mc {\hat Y}$. 
The nodes corresponding to partial assignments in $\mc {\hat Y}$ are required to be further explored to check whether or not the value of $w^*$ can be improved. 
That is, the above process is repeated for every node $v$ with partial assignment in $\mc {\hat Y}$ as the $r$-$v$ path is fixed in the new restricted/relaxed DDs. 
The algorithm terminates when $\mc {\hat Y}$ becomes empty, at which point $w^*$ is the optimal value.       

\section{DD-BD Formulation for the SGUFP} \label{sec: DDBD}
  
In this section, we adapt the DD-BD framework described in Section~\ref{subsec: Continuous DD Models} to solve the SGUFP.
  
\subsection{MIP Formulation} \label{subsec: MIP Formulation}

We study the MIP formulation of the SGUFP based on that of its deterministic counterpart given in \cite{davarnia2019network}.
Consider a network $G=(V,A)$ with node set $V \coloneqq V'\cup\{s,t\}$ and arc set $A$, where $s$ and $t$ are source and sink nodes, respectively. 
The source node is connected to all the supply nodes in $S \subseteq V'$, and the sink node is connected to all the demand nodes in $D \subseteq V'$.
Figure~\ref{fig:network} illustrates the general structure of this network.
For a node $q\in V$, let $\delta^-(q) \coloneqq \{i\in V \suchthat (i,q)\in A\}$ and $\delta^+(q) \coloneqq \{j\in V \suchthat (q,j)\in A\}$ show the set of incoming and outgoing neighbors of $q$, respectively. 
Define $\bar V \subseteq V'$ as a subset of vertices that must satisfy the NSNM requirement. 
For each node $q \in \bar V$, let binary variable $y^q_{ij}\in \{0,1\}$ represent whether or not the flow entering node $q\in \bar V$ through arc $(i,q)$ leaves node $q$ through arc $(q,j)$.
The first stage of SGUFP determines the matching pairs between incoming and outgoing arcs of unsplittable nodes as follows:

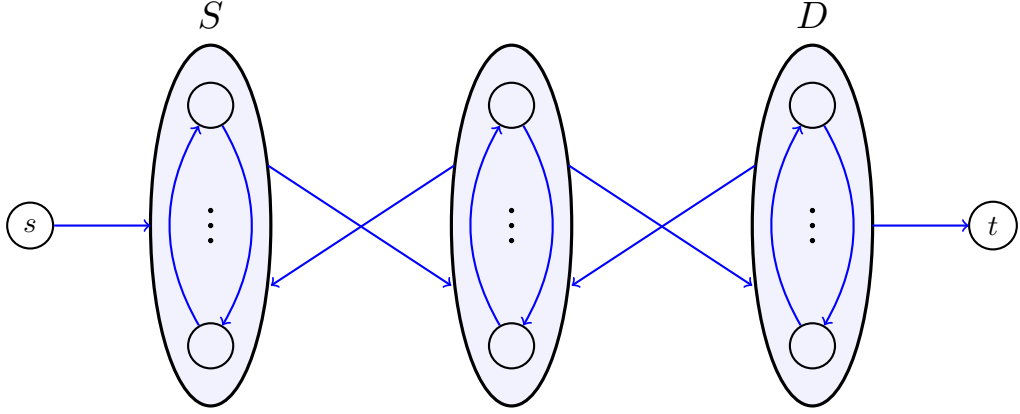
\begin{figure}
\centering
\begin{tikzpicture}[scale=0.8]
\begin{scope}[every node/.style={circle,thick,draw,minimum size=0.6cm}]
\filldraw[fill=blue!5, very thick](-1,0) ellipse (1 and 3);
\filldraw[fill=blue!5, very thick](4,0) ellipse (1 and 3);
\filldraw[fill=blue!5, very thick](9,0) ellipse (1 and 3);

\node[draw=none]() at (-1,3.5) {\large $S$};
\node[draw=none]() at (9,3.5) {\large $D$};

\node(s) at (-4,0) {$s$};
\node(t) at (12,0) {$t$};

\node(s1) at (-1,2) {};
\filldraw[black] (-1,0.25) circle (1pt);
\filldraw[black] (-1,0) circle (1pt);
\filldraw[black] (-1,-0.25) circle (1pt);
\node(s2) at (-1,-2) {};

\node [label=above:{}] (n1) at (4,2) {};
\filldraw[black] (4,0.25) circle (1pt);
\filldraw[black] (4,0) circle (1pt);
\filldraw[black] (4,-0.25) circle (1pt);
\node [label=above:{}] (n2) at (4,-2) {};

\node [label={[label distance=-0.3cm]}] (d1) at (9,2) {};
\filldraw[black] (9,0.25) circle (1pt);
\filldraw[black] (9,0) circle (1pt);
\filldraw[black] (9,-0.25) circle (1pt);
\node [label={[label distance=-0.3cm]}] (d2) at (9,-2) {};

\path [->, thick, color=blue] (s1) edge[bend left=30] (s2);
\path [->, thick, color=blue] (s2) edge[bend left=30] (s1);
\path [->, thick, color=blue] (n1) edge[bend left=30] (n2);
\path [->, thick, color=blue] (n2) edge[bend left=30] (n1);
\path [->, thick, color=blue] (d1) edge[bend left=30] (d2);
\path [->, thick, color=blue] (d2) edge[bend left=30] (d1);
\draw[->, thick, color=blue] (-3.6,0) -- (-2,0);
\draw[->, thick, color=blue] (-0.05,1) -- (3,-1);
\draw[<-, thick, color=blue] (0,-1) -- (3.05,1);
\draw[->, thick, color=blue] (4.95,1) -- (8,-1);
\draw[<-, thick, color=blue] (5,-1) -- (8.05,1);
\draw[->, thick, color=blue] (10,0) -- (11.6,0);
\end{scope}
\end{tikzpicture}
\caption{Illustration of network $G=(V' \cup \{s,t\},A)$}
\label{fig:network}
\end{figure}

\begin{subequations}
\begin{align}
\max \quad &z \label{obj1} \\  
\text{s.t.} \quad &\sum_{j\in \delta^+(q)}y^q_{ij} \le 1 &\forall i\in \delta^-(q),~\forall q\in \bar V \label{single_match_1} \\
&\sum_{i\in \delta^-(q)}y^q_{ij} \le 1 &\forall j\in \delta^+(q),~\forall q\in \bar V \label{single_match_2} \\
&y^q_{ij}\in \{0,1\} &\forall (i,j)\in \delta^-(q) \times \delta^+(q),~\forall q\in \bar V, \label{binary_y}
\end{align}
\end{subequations}

where constraints~\eqref{single_match_1} ensure that each incoming arc to a node with NSNM requirement is assigned to at most one outgoing arc, and constraints~\eqref{single_match_2} guarantee that each outgoing arc from such a node is matched with at most one incoming arc.

In \eqref{obj1}--\eqref{binary_y}, variable $z$ represents the objective value of the second stage of SGUFP where the demand uncertainty is taken into account.
This demand uncertainty is modeled by a set $\Xi$ of scenarios for the demand vector $\vc{d}^{\xi}$ with occurrence probability $\prob^\xi$ for each scenario $\xi \in \Xi$.
Let continuous variable $x^\xi_{ij}\in \mathbb{R}_+$ denote the flow from node $i$ to node $j$ through arc $(i,j)$ under scenario $\xi \in \Xi$.
We further assign a reward $r_{ij}$ per unit flow to be collected by routing flow through arc $(i,j)$.
It follows that $z = \sum_{\xi \in \Xi}\prob^\xi z^{\xi}$, where $z^{\xi}$ is the objective value of the second stage of SGUFP for each scenario $\xi \in \Xi$.
This subproblem is formulated as follows for a given $\vc{y}$ vector:

\begin{subequations}
\begin{align}
\max \quad &\sum_{q\in V}\sum_{j\in \delta^+(q)}r_{qj}x^\xi_{qj} \label{obj} \\  
\text{s.t.} \quad &\sum_{i\in \delta^-(q)}x^\xi_{iq} - \sum_{j\in \delta^+(q)}x^\xi_{qj} = 0 &\forall q\in V' \label{flow_conservation} \\
&\ell^\xi_{iq} \le x^\xi_{iq} \le u^\xi_{iq} &\forall i\in \delta^-(q),~\forall q\in V \label{flow_bounds} \\
&x^\xi_{iq}-x^\xi_{qj} \le u^\xi_{iq}(1-y^q_{ij}) &\forall(i,j)\in \delta^-(q)\times \delta^+(q),~\forall q\in{\bar V} \label{NSNM_1} \\
&x^\xi_{qj}-x^\xi_{iq} \le u^\xi_{qj}(1-y^q_{ij}) &\forall(i,j)\in \delta^-(q) \times \delta^+(q),~\forall q\in \bar V \label{NSNM_2} \\
&x^\xi_{iq} \le u^\xi_{iq} \sum_{j\in \delta^+(q)}y^q_{ij} &\forall i\in \delta^-(q),~\forall q\in \bar V \label{NSNM_3} \\
&x^\xi_{qj} \le u^\xi_{qj} \sum_{i\in \delta^-(q)}y^q_{ij} &\forall j\in \delta^+(q),~\forall q\in \bar V \label{NSNM_4} \\
&x^\xi_{ij} \geq 0 &\forall (i,j)\in A. \label{x_bound}
\end{align}
\end{subequations}

In the above formulation, the objective function captures the total reward collected by routing flows throughout the network (from the source $s$ to the sink $t$) to satisfy demands.  
The flow-balance requirements are represented by \eqref{flow_conservation}. 
Constraints~\eqref{flow_bounds} bound the flow on each arc from below and above. 
To impose the demand requirement for each scenario $\xi \in \Xi$, we fix $\ell^{\xi}_{qt}=u^{\xi}_{qt}=d^{\xi}_q$ for all demand nodes $q\in D$ with demand $d^{\xi}_q$, and leave the lower and upper bound values unchanged for all other arcs. Constraints~\eqref{NSNM_1}--\eqref{NSNM_4} model the NSNM requirement for each node $q \in \bar{V}$.
In particular, \eqref{NSNM_1} and~\eqref{NSNM_2} ensure that matching arcs $(i,q)$ and $(q,j)$ have equal flows. 
Constraints~\eqref{NSNM_3} and~\eqref{NSNM_4} guarantee that an arc without a matching pair does not carry any flow.  
We note here that the Constraint \eqref{flow_conservation} is implied by other constraints of the above subproblem under the assumption that $\vc{y}$ is feasible to the master problem \eqref{obj1}--\eqref{binary_y}.
However, we maintain this constraint in the subproblem because the master formulation in our DD-based approach, as will be described in Section~\ref{subsec:master}, may produce a solution that is not feasible to \eqref{obj1}--\eqref{binary_y}.
As a result, the addition of the Constraint \eqref{flow_conservation} will lead to a tighter subproblem formulation.

As discussed in Section~\ref{subsec: Continuous DD Models}, the first step to use the DD-BD algorithm is to decompose the underlying problem into a master and a subproblem. 
The above two-stage formulation of the SGUFP is readily amenable to BD since the first stage problem \eqref{obj1}-\eqref{binary_y} can be considered as the master problem together with some valid lower and upper bounds $-\Gamma$ and $\Gamma$ on $z$ induced from the boundedness of the MIP formulation. 
For a given $\vc{y}$ value obtained from the master problem and a scenario $\xi \in \Xi$, the second stage problem \eqref{obj}-\eqref{x_bound} can be viewed as the desired subproblems.
The optimality/feasibility cuts obtained from each scenario-based subproblem are then added to the master problem through aggregation as described in Section~\ref{sec:subFormulation}.

\subsection{DD-BD: Master Problem Formulation} \label{subsec:master}
While the DD-BD Algorithm~\ref{alg: DD-BD} provides a general solution framework for any bounded MIP, its DD component is problem-specific, i.e., it should be carefully designed based on the specific structure of the underlying problem.
In this section, we design such an oracle for the SGUFP that represents the feasible region $\{\eqref{single_match_1}-\eqref{binary_y}, z\in[-\Gamma, \Gamma]\}$ of the master problem~\eqref{obj1}-\eqref{binary_y}. 
To model this feasible region in the original space of $(\vc y;z)$ variables, a DD would require $\sum_{q \in \bar V}|\delta^-(q)| \times |\delta^+(q)|$ arc layers to represent binary variables $\vc y$ and one arc layer to encode the continuous variable $z$. 
Constructing such a DD, however, would be computationally cumbersome due to the large number of the arc layers.
To mitigate this difficulty, we take advantage of the structural flexibility of DDs in representing \textit{irregular} variable types that cannot be used in standard MIP models.
One such variable type is the index set, where arc layers represent indices, rather than domain values. 
We next show that we can remarkably reduce the number of DD arc layers by reformulating the master problem in a transformed space of variables defined over index sets. 

Consider a node $q\in \bar V$.
In the following, we define mappings that assign an index to each incoming and outgoing arc of $q$. These mappings enable us to define new variables to reduce the number of DD arc layers.
Let $\ind^-(i,q)$ be a one-to-one mapping from incoming arcs $(i,q)$, for $i \in \delta^-(q)$, to the index set $\{1,2,\dots,|\delta^-(q)|\}$.
Similarly, let $\ind^+(q,j)$ be a one-to-one mapping from outgoing arcs $(q,j)$, for $j \in \delta^+(q)$, to the index set $\{1,2,\dots,|\delta^+(q)|\}$.
For each incoming arc $(i,q)$ with index $h = \ind^-(i,q)$, we define an integer variable $w^q_h \in \{0,1,\dots,|\delta^+(q)|\}$ such that $w^q_h = 0$ if this incoming arc is not paired with any outgoing arc, and $w^q_h = k > 0$ if this arc is matched with an outgoing arc $(q,j)$ with index $k = \ind^+(q,j)$. 

Next, we give a formulation in the space of $\vc{w}$ variables that describes the matching between incoming and outgoing arcs of $q$ for all $q \in \bar V$.
In the following, $\sign(.)$ represents the sign function that returns $1$ if its argument is strictly positive, $0$ if the argument is zero, and $-1$ otherwise. 
Further, the operator $|.|$, when applied on a set, represents the set size; and when applied on a real number, it represents the absolute value.

\begin{proposition} \label{prop: w}
Formulation
\begin{subequations}
\begin{align}
&\sum_{i \in \delta^-(q)} \sign\left(\left|w^q_{\ind^-(i,q)} - \ind^+(q,j)\right|\right) \geq \left|\delta^-(q)\right| - 1 &\forall j \in \delta^+(q),~\forall q\in \bar V \label{eq:w1}\\
& w^q_{\ind^-(i,q)} \in \left\{0,1,\dots,\left|\delta^+(q)\right|\right\} & \forall i \in \delta^-(q), ~\forall q\in \bar V \label{eq:w2}
\end{align}
\end{subequations}
models the matching between incoming and outgoing arcs of nodes $q \in \bar{V}$.
\end{proposition}

\proof{Proof.}
We show the result for a single node $q \in \bar{V}$.
The extension to the multiple node case is straightforward as the matching problem for each node is independent from other nodes. 
For the direct implication, assume that $M^q$ is a matching between incoming and outgoing arcs of $q$, with elements of the form $(i,j)$ that represent a matching between the incoming arc $(i,q)$ and the outgoing arc $(q,j)$.
We show that variables $\vc{w}$ associated with matching pairs in $M^q$ satisfy constraints \eqref{eq:w1} and \eqref{eq:w2}.
It follows from the definition of $\vc{w}$ that, for each $(i,j) \in M^q$, we have $w^q_{\ind^-(i,q)} = \ind^+(q,j)$.
Also, for any $i \in \delta^-(q)$ that does not have a matching pair in $M^q$, we have $w^q_{\ind^-(i,q)} = 0$.
These value assignments show that $\vc{w}$ satisfies \eqref{eq:w2} as the image of $\ind^+$ mapping is $\{1,\dots,|\delta^+(q)|\}$.  
For each $i \in \delta^-(q)$ and $j \in \delta^+(q)$, we have $\left|w^q_{\ind^-(i,q)} - \ind^+(q,j)\right| \geq 0$, with equality holding when $(i,j) \in M^q$. 
For each $j \in \delta^+(q)$, there are two cases.
For the first case, assume that $(i,j) \notin M^q$ for any $i \in \delta^-(q)$.
As a result, $\left|w^q_{\ind^-(i,q)} - \ind^+(q,j)\right| > 0$ for all $i \in \delta^-(q)$.
Applying the $\sign(.)$ function on these terms yields $\sign\left(\left|w^q_{\ind^-(i,q)} - \ind^+(q,j)\right|\right) = 1$, which implies that
$\sum_{i \in \delta^-(q)} \sign\left(\left|w^q_{\ind^-(i,q)} - \ind^+(q,j)\right|\right) = \left|\delta^-(q)\right|$, satisfying \eqref{eq:w1}.
For the second case, assume that $(i^*,j) \in M^q$ for some $i^* \in \delta^-(q)$.
As a result, we have
$\sum_{i \in \delta^-(q)} \sign\left(\left|w^q_{\ind^-(i,q)} - \ind^+(q,j)\right|\right) = \left|\delta^-(q)\right| - 1$ since
$\sign\left(\left|w^q_{\ind^-(i^*,q)} - \ind^+(q,j)\right|\right) = \left|w^q_{\ind^-(i^*,q)} - \ind^+(q,j)\right| = 0$, satisfying \eqref{eq:w1}.

For the reverse implication, assume that $\vc{w}$ is a feasible solution to \eqref{eq:w1}--\eqref{eq:w2}. We show that the pairs of the form $(i,j)$ encoded by these variables constitute a feasible matching between incoming and outgoing arcs of $q$, i.e., (i) each arc $(i,q)$ is matched with at most one arc $(q, j)$, and (ii) each arc $(q,j)$ is matched with at most one arc $(i,q)$.
It follows from constraint \eqref{eq:w2} that, for each $i \in \delta^-(q)$, variable $w^q_{\ind^-(i,q)}$ takes a value between $\{0,1,\dots,|\delta^+(q)|\}$.
If $w^q_{\ind^-(i,q)} = 0$, then $(i,q)$ is not matched with any outgoing arc, otherwise it is matched with arc $(q,j)$ with $\ind^+(q,j) = w^q_{\ind^-(i,q)}$.
This ensures that condition (i) above is satisfied for this matching collection.
Further, for each $j \in \delta^-(q)$, constraint \eqref{eq:w1} implies that $\sign\left(\left|w^q_{\ind^-(i,q)} - \ind^+(q,j)\right|\right)$ can be equal to zero for at most one $i \in \delta^-(q)$.
In such a case, we would have at most one matching pair of the form $(i,j)$ in the collection, showing that condition (ii) above is satisfied.
\Halmos
\endproof

It follows from Proposition~\ref{prop: w} that constraints \eqref{eq:w1}-\eqref{eq:w2} can replace \eqref{single_match_1}-\eqref{binary_y} in the master problem \eqref{obj1}-\eqref{binary_y} to obtain the following master problem in a transformed space of variables.
\begin{align}
\max_{\vc w; z}\left\{z \suchthat \eqref{eq:w1}-\eqref{eq:w2}, z\in[-\Gamma, \Gamma]\right\}. \label{eq:master}
\end{align}

Note that formulation \eqref{eq:master} is an integer nonlinear program (INLP) with nonconvex and noncontinuous constraint functions.
Such a formulation is extremely difficult for conventional MINLP techniques and solvers to handle.
However, due to structural flexibility of DDs in representing integer nonlinear programs, this problem can be easily modeled via a DD; see \cite{davarnia2020outer} for a detailed account on using DDs for modeling INLPs. 
In the following, we present an algorithm to construct DDs in the space of $(\vc w;z)$ variables for the master problem \eqref{eq:master} with a single node $q \in \bar V$. 
The extension to the case with multiple nodes follows by replicating the DD structure. 
The output of Algorithm~\ref{alg: DD master} is a DD with $|\delta^-(q)|+1$ arc layers where the first $|\delta^-(q)|$ layers represent $\vc w$ variables and the last layer encodes variable $z$.
In this algorithm, $s_u$ denotes the state value of DD node $u$.
The core idea of the algorithm is to use unpaired outgoing arcs of $q$ as the state value at each DD layer that represents the matching for an incoming arc of $q$.

\begin{algorithm}[!ht] 
\caption{Construction of DD for the master problem of SGUFP with a node $q\in \bar V$} \label{alg: DD master}
\KwData{node $q \in \bar V$, parameter $\Gamma$ }
\KwResult{an exact DD $\mc D$}

create the root node $r \in \mc{U}_1$ with state $s_r = \{0,1,\dots,|\delta^+(q)|\}$

\ForAll{$i \in\{1,2,\dots,|\delta^-(q)|\}$ and $u \in \mc U_i$}
{
\ForAll{$\ell \in s_u$}{create a node $v\in \mc U_{i+1}$ with state $(s_u \setminus \{\ell\})\cup \{0\}$ and an arc $a \in \mc A_i$ connecting $u$ to $v$ with label $\l(a)=\ell$} 
}
\ForAll{$u\in \mc U_{1+|\delta^-(q)|}$}{ create two arcs $a_1, a_2 \in \mc A_{1+|\delta^-(q)|}$ connecting $u$ to the terminal node with labels $l(a_1) = \Gamma$ and $l(a_2) = -\Gamma$.}
\end{algorithm}

Next, We show that the solution set of the DD constructed by Algorithm~\ref{alg: DD master} \textit{represents} the feasible region of $\eqref{eq:master}$.
Note here that DD representation of a MIP set, as described in Section~\ref{subsec: Continuous DD Models}, does not imply the encoding of all of the solutions of the set, but rather the encoding of a subset of all solutions that subsumes all the extreme points of the set. 
Such a representation is sufficient to solve an optimization problem over the set with an objective function convex in continuous variables, which is the case for $\eqref{eq:master}$.

\begin{theorem} \label{thm:master exact}
Consider a SGUFP with $\bar{V} = \{q\}$.
Let $\mc D$ be a DD constructed by Algorithm~\ref{alg: DD master}. 
Then, $\Sol(\mc{D})$ represents the feasible region of \eqref{eq:master}.
\end{theorem}

\proof{Proof.}
$(\subseteq)$ Consider an $r$-$t$ path of $\mc D$ that encodes solution $(\tilde{\vc w}^q,z)$. 
According to Algorithm~\ref{alg: DD master}, the labels of the first $|\delta^-(q)|$ arcs of this path belong to $\{0,1,\dots,|\delta^+(q)|\}$, showing that $\tilde{\vc{w}}^q$ satisfies constraints $\eqref{eq:w2}$.
Assume by contradiction that $\tilde{\vc w}^q$ does not satisfy constraints~\eqref{eq:w1}, i.e., $\sum_{i \in \delta^-(q)} \sign\left(\left|w^q_{\ind^-(i,q)} - \ind^+(q,j)\right|\right) \leq \left|\delta^-(q)\right| - 2$ for some $j \in \delta^+(q)$.
This implies that $\tilde{w}^q_{\ind^-(i',q)} = \tilde{w}^q_{\ind^-(i'',q)}  = \ind^+(q,j)$ for two distinct $i', i'' \in \delta^-(q)$.
In other words, the arcs at layers $\ind^-(i',q)$ and $\ind^-(i'',q)$ of the selected $r$-$t$ path both share the same label value $\ind^+(q,j)$.
According to line 3 of Algorithm~\ref{alg: DD master}, we must have that the state value of nodes at layers $\ind^-(i',q)$ and $\ind^-(i'',q)$ of the $r$-$t$ path both contain $\ind^+(q,j)$.
This is a contradiction to the state update policy in line 4 of Algorithm~\ref{alg: DD master}, since positive arc labels at each layer of the DD will be excluded from the state value of the subsequent nodes.

$(\supseteq)$ Consider a feasible solution point $(\tilde{\vc w}^q; \tilde z)$ of \eqref{eq:master}.
Suppose $\tilde{\vc w}^q = (\ell_1,\ell_2,\dots,\ell_{|\delta^-(q)|})$. According to constraints~\eqref{eq:w1}, no two coordinates of $\tilde{\vc w}^q$ have the same positive value. 
The state value at the root node in $\mc{D}$ contains all index values $\{0,1,\dots,|\delta^+(q)|\}$.
According to Algorithm~\ref{alg: DD master}, there exists an arc with label $\ell_1$ at the first layer of $\mc{D}$.
The state value at the head node of this arc, therefore, contains $\ell_2 \in \{0,1,\dots,|\delta^+(q)|\} \setminus \{\ell_1\}$, which guarantees an arc with label $\ell_2$ at the second layer of this path.
Following a similar approach, we can track a path from the root to layer $|\delta^-(q)|$ whose arcs labels match values of $\tilde{\vc{w}}^q$.
Note for the last layer that $\tilde z \in [-\Gamma, \Gamma]$, which is included in the interval between arc labels of the last layer of $\mc{D}$.
As a result, $(\tilde{\vc w}^q; \tilde z)$ is represented by an $r$-$t$ path of $\mc{D}$.
\Halmos
\endproof

The main purpose of using a DD that models the master problem \eqref{eq:master} over one that models~\eqref{obj1}-\eqref{binary_y} is the size reduction in arc layers that represent variables $\vc{w}$ as compared with variables $\vc{y}$.
It turns out that this space transformation can significantly improve the solution time of the DD approach.
We refer the interested reader to Appendix~\ref{app:comparison} for a detailed discussion on these advantages, including preliminary computational results.

Constructing exact DDs as described in Algorithm~\ref{alg: DD master} can be computationally expensive for large size problems.
As discussed in Section~\ref{subsec: Continuous DD Models}, relaxed and restricted DDs are used to circumvent this difficulty.
Building restricted DDs is straightforward as it involves the selection of a subset of $r$-$t$ paths of the exact DD that satisfy a preset width limit.
Constructing relaxed DDs, on the other hand, requires careful manipulation of the DD structure to merge nodes in such a way that it encodes a superset of all $r$-$t$ paths of the exact DD.
We demonstrate a method to construct such relaxed DDs in Algorithm~\ref{alg: relaxed DD}.
Similarly to Algorithm~\ref{alg: DD master}, this algorithm is presented for a single NSNM node, but can be extended to multiple nodes by replicating the procedure.

\begin{algorithm}[!ht] 
\caption{Construction of relaxed DD for the master problem of SGUFP with a node $q\in \bar V$} \label{alg: relaxed DD}
\KwData{node $q \in \bar V$, parameter $\Gamma$ }
\KwResult{a relaxed DD $\overline{\mc D}$}

create the root node $r \in \mc{U}_1$ with state $s_r = \{0,1,\dots,|\delta^+(q)|\}$

\ForAll{$i \in\{1,2,\dots,|\delta^-(q)|\}$ and $u \in \mc U_i$}
{
\ForAll{$\ell \in s_u$}{create a node $v\in \mc U_{i+1}$ with state $(s_u \setminus \{\ell\})\cup \{0\}$ and an arc $a \in \mc A_i$ connecting $u$ to $v$ with label $\l(a)=\ell$} 

select a subset of nodes $v_1,v_2,\dots,v_k \in \mc{U}_{i+1}$ and merge them into node $v'$ with state $s_{v'} = \bigcup_{j=1}^{k}s_{v_j}$ 

}
 
\ForAll{$u\in \mc{U}_{1+|\delta^-(q)|}$}{create two arcs $a_1, a_2 \in \mc{A}_{1+|\delta^-(q)|}$ connecting $u$ to the terminal node with labels $l(a_1) = \Gamma$ and $l(a_2) = -\Gamma$.}

\end{algorithm}

\begin{theorem}
Consider a SGUFP with $\bar V = \{q\}$.
Let $\overline{\mc{D}}$ be a DD constructed by Algorithm~\ref{alg: relaxed DD}.
Then, $\overline{\mc{D}}$ represents a relaxation of the feasible region of~\eqref{eq:master}.
\end{theorem}

\proof{Proof.}
Let $\dot{\mc D}$ be the DD constructed by Algorithm~\ref{alg: DD master} for the master problem~\eqref{eq:master} with a single node $q \in \bar V$. 
It suffices to show that the solution set of $\overline{\mc D}$ provides a relaxation for that of $\dot{\mc D}$. 
Pick a root-terminal path $\dot P$ of $\dot{\mc D}$ with encoding point $(\dot{\vc w}^q;\dot z)$. 
We show that there exist a root-terminal path $\overline{P}$ of $\overline{\mc D}$ with encoding point $(\overline{\vc w}^q;\overline{z})$ such that $\overline{\vc w}^q = \dot{\vc w}^q$ and $\overline{z} = \dot{z}$.
Given a DD, define $P_k$ to be a sub-path composed of arcs in the first $k$ layers, for $1 \leq k \le |\delta^-(q)|$. 
We show for any sub-path $\dot{P}_k$ of $\dot{\mc D}$ with encoding point $\dot{\vc w}^q_k = (\dot w^q_1,\dots, \dot w^q_k)$, there exists a sub-path $\overline{P}_k$ of $\overline{\mc D}$ with encoding point $\overline{\vc w}_k = (\overline{w}_1, \dots, \overline{w}_k)$ such that $\overline{\vc w}_h = \dot{\vc w}_h$ for $h = 1, \dots, k$. 
Note that we only need to prove the matching values for $k \le |\delta^-(q)|$, because each node at node layer $|\delta^-(q)|+1$ of both $\dot{\mc D}$ and $\overline{\mc D}$ is connected by two arcs with labels $-\Gamma$ and $\Gamma$ to the terminal node, and thus there are always matching arcs with the same label for the last layer, i.e., $\overline{z} = \dot{z}$. 
We prove the result by induction on $k$. 
The base case for $k=1$ is trivial, since $\overline{\mc D}$ contains arcs with labels $\{0,1,\dots,|\delta^+(q)|\}$ in the first layer, which includes the label value of the first arc on $\dot{P}_1$. 
For the induction hypothesis, assume that the statement is true for $k=d$, i.e., for the sub-path $\dot{P}_d$ with label values $\dot{\vc{w}}^q_d = (\dot w^q_1,\dots, \dot w^q_d)$, there is sub-path $\overline{P}_d$ of $\overline{\mc D}$ with matching arc labels.
We show the statement holds for $d+1$. 
Let $u \in \dot{A}_{d+1}$ and $v \in \overline{A}_{d+1}$ be the end nodes of $\dot{P}_d$ and $\overline{P}_d$, respectively. 
It follows from Algorithm~\ref{alg: DD master} that the index set representing the state value at node $u$ contains $\dot w^q_{d+1}$, i.e., $\dot w^q_{d+1} \in \dot s_u = \{0\} \cup \{1, \dots, |\delta^+(q)|\} \setminus \{\dot{w}_1, \dot{w}_2, \dots, \dot{w}_d \}$. 
The merging step in line 5 of Algorithm~\ref{alg: relaxed DD}, on the other hand, implies that $\overline{s}_v \supseteq \{0\} \cup \{1, \dots, |\delta^+(q)|\} \setminus \{\overline{w}_1, \overline{w}_2, \dots, \overline{w}_d \} = \{0\} \cup \{1, \dots, |\delta^+(q)|\} \setminus \{\dot{w}_1, \dot{w}_2, \dots, \dot{w}_d \} = \dot s_u$, where the inclusion follows from the fact that state values at nodes on path $\overline P_d$ contain those of each individual path due to merging operation, and the first equality holds because of the induction hypothesis.
As a result, $\overline{s}_v$ must contain $\dot w^q_{d+1}$, which implies that there exists an arc with $\dot w^q_{d+1}$ connected to node $v$ on $\overline P_d$.
Attaching this arc to $\overline P_d$, we obtain the desired sub-path $\overline P_{d+1}$.
\Halmos
\endproof

\subsection{DD-BD: Subproblem Formulation} \label{sec:subFormulation}

At each iteration of the DD-BD algorithm, an optimal solution of the master problem is plugged into the subproblems to obtain feasibility/optimality cuts.
For the SGUFP formulation, this procedure translates to obtaining an optimal solution of \eqref{eq:master} in the space of $\vc w$ variables, which is used to solve the subproblem~\eqref{obj}-\eqref{x_bound}.
The formulation of the subproblem, however, is defined over the original binary variables $\vc y$, and the resulting feasibility/optimality cuts are generated in this space.
To remedy this discrepancy between the space of variables in the master and subproblems, we need to find a one-to-one mapping between variables $\vc w$ and $\vc y$, as outlined next.

\begin{proposition} \label{prop:mapping}
Consider a node $q \in \bar V$.
Let $\vc{y}^q$ be a feasible solution to \eqref{single_match_1}-\eqref{binary_y}.
Then, $\vc{w}^q$ obtained as 
\begin{align}
&w^q_{\ind^-(i,q)} = \sum_{j \in \delta^+(q)}  \ind^+(q,j) \, y^q_{ij} &\forall i \in \delta^-(q), \label{eq:y-w}
\end{align}
is a feasible solution to \eqref{eq:w1}-\eqref{eq:w2}.
Conversely, let $\vc{w}^q$ be a feasible solution to \eqref{eq:w1}-\eqref{eq:w2}.
Then, $\vc{y}^q$ obtained as 
\begin{align}
&y^q_{ij} = 1 - \sign\left(\left|w^q_{\ind^-(i,q)} - \ind^+(q,j) \right|\right) &\forall (i,j) \in \delta^-(q) \times \delta^+(q), \label{eq:w-y}
\end{align}
is a feasible solution to \eqref{single_match_1}-\eqref{binary_y}.
\end{proposition}

\proof{Proof.}
For the direct statement, let $\vc y^q$ be a feasible solution to~\eqref{single_match_1}-\eqref{binary_y}, and construct a vector $\vc w^q$ according to~\eqref{eq:y-w}. 
We show that $\vc w^q$ satisfies all constraints~\eqref{eq:w1}-\eqref{eq:w2}. 
First, we show that constraints~\eqref{eq:w1} are satisfied. 
Assume by contradiction that there exists $j' \in \delta^+(q)$ such that $\sum_{i \in \delta^-(q)} \sign\left(\left|w^q_{\ind^-(i,q)} - \ind^+(q,j')\right|\right) \le |\delta^-(q)| - 2$. 
This implies that $w^q_{\ind^-(i',q)} = w^q_{\ind^-(i'',q)} = \ind^+(q,j')$ for some $i',i'' \in \delta^-(q)$.
Then, we can write that
\begin{align*}
w^q_{\ind^-(i',q)} = \sum_{j \in \delta^+(q)} \ind^+(q,j)y^q_{i'j} = \ind^+(q,j') = \sum_{j \in \delta^+(q)} \ind^+(q,j)y^q_{i''j} = w^q_{\ind^-(i'',q)},  
\end{align*}
where the first and last equalities hold by~\eqref{eq:y-w}. 
The second and third equalities in the above chain of relations imply that $y^q_{i'j'} =y^q_{i''j'}=1$, since $\ind^+(q,j') > 0$. 
This violates constraints~\eqref{single_match_2}, reaching a contradiction.
Next, we show that constraints~\eqref{eq:w2} are satisfied.
The proof follows directly from construction of $\vc{w}^q$ and constraints~\eqref{single_match_1}.       

For the converse statement, let $\vc w^q$ be a feasible solution to~\eqref{eq:w1}-\eqref{eq:w2}, and construct a vector $\vc y^q$ according to~\eqref{eq:w-y}. 
We show that $\vc y^q$ satisfies all constraints~\eqref{single_match_1}-\eqref{binary_y}.
To show that each constraint~\eqref{single_match_1} is satisfied, consider $i \in \delta^-(q)$.
We can write that
\begin{align*}
\sum_{j \in \delta^+(q)} y^q_{ij} = |\delta^+(q)| - \sum_{j \in \delta^+(q)} \sign\left(\left|w^q_{\ind^-(i,q)} - \ind^+(q,j)\right|\right) \le |\delta^+(q)| - \left(|\delta^+(q)| - 1\right) = 1,   
\end{align*}
where the first equality follows from the construction of $\vc{y}^q$, and the inequality holds by~\eqref{eq:w2} as $\left|w^q_{\ind^-(i,q)} - \ind^+(q,j)\right| = 0$ for at most one index $j \in \delta^+(q)$. 
To show that each constraint~\eqref{single_match_2} is satisfied, select $j \in \delta^+(q)$. 
We have
\begin{align*}
\sum_{i \in \delta^-(q)} y^q_{ij} = |\delta^-(q)| - \sum_{i \in \delta^-(q)} \sign\left(\left|w^q_{\ind^-(i,q)} - \ind^+(q,j)\right|\right) \le 1, 
\end{align*}
where the equality follows from the construction of $\vc{y}^q$, and the inequality holds because of constraint~\eqref{eq:w1}.
Finally, each constraint~\eqref{binary_y} is satisfied due to the fact that $1- \sign(|.|) \in \{0,1\}$. 
\Halmos
\endproof

\begin{proposition} \label{prop:1-1}
Mappings described by~\eqref{eq:y-w} and~\eqref{eq:w-y} are one-to-one over their respective domains.
\end{proposition}

\proof{Proof.}
Note that the mapping described by \eqref{eq:y-w} is a linear transformation of the form $\vc{w}^q = B \vc{y}^q$ with coefficient matrix $B \in \mathbb{Z}^{|\delta^-(q)| \times ( |\delta^-(q)| |\delta^+(q)|)}$. 
It is clear from the identity block structure of $B$, that it is full row-rank, since each column contains a single non-zero element while each row has at least one non-zero element.
As a result, the null space of $B$ is the origin, which implies that $\hat{\vc{w}}^q = \tilde{\vc{w}}^q$ only if $\hat{\vc{y}}^q = \tilde{\vc{y}}^q$.

For the mapping described by \eqref{eq:w-y}, let distinct points $\hat{\vc w}^q$ and $\tilde{\vc w}^q$ satisfy \eqref{eq:w2}.
Construct vectors $\hat{\vc y}^q$ and $\tilde{\vc y}^q$ by~\eqref{eq:w-y} using $\hat{\vc w}^q$ and $\tilde{\vc w}^q$, respectively. 
Because $\hat{\vc w}^q$ and $\tilde{\vc w}^q$ are distinct, there must exist $i \in \delta^-(q)$ such that $\hat w^q_{\ind^-(i,q)} \neq \tilde w^q_{\ind^-(i,q)}$.
This implies that at least one of these variables, say $\hat w^q_{\ind^-(i,q)}$, is non-zero.
It follows from \eqref{eq:w2} that $\hat w^q_{\ind^-(i,q)} = \ind^+(q,j')$ for some $j' \in \delta^+(q)$, and that $\hat w^q_{\ind^-(i,q)} \neq \ind^+(q,j')$.
According to \eqref{eq:w-y}, we write that $\hat y_{ij'} = 1 - \sign\left(\left|\hat w^q_{\ind^-(i,q)} - \ind^+(q,j') \right|\right) = 1$, and that $\tilde y_{ij'} = 1 - \sign\left(\left|\tilde w^q_{\ind^-(i,q)} - \ind^+(q,j') \right|\right) = 0$, showing that $\hat{\vc y}^q \neq \tilde{\vc y}^q$.
\Halmos
\endproof

Using the results of Propositions~\ref{prop:mapping} and \ref{prop:1-1}, we can apply the DD-BD Algorithm~\ref{alg: DD-BD} in its entirety for the SGUFP. 
In particular, at each iteration of the algorithm, we can transform the optimal solution $(\bar{\vc{w}},\bar{z})$ obtained from the DD representing the master problem \eqref{eq:master} into a solution $(\bar{\vc{y}},\bar{z})$ through the mapping~\eqref{eq:w-y}.
Given an optimal first-stage solution $\bar{\vc y}$, we can solve $|\Xi|$ separate subproblems; one for each demand realization in the second-stage. The feasibility cuts obtained from subproblems, which are in the space of $\vc y$ variables, are translated back into the space of $\vc w$ variables through the mapping~\eqref{eq:y-w} and added to the master problem. 
Further, in a case where all subproblems produce an optimality cut, they can be aggregated to generate an optimality cut in the space of $(\vc y, z)$, which is added to the master problem after being translated into the space of $(\vc w,z)$ variables.
The master DD will be refined with respect to the resulting inequalities, and an optimal solution is returned to be used for the next iteration.

In the remainder of this section, we present details on the derivation of optimality/feasibility cuts from subproblem \eqref{obj}-\eqref{x_bound}. 
Consider the following partitioning of the set of arcs $A$ into subsets     
\begin{align*}
&A_1 \coloneqq \left\{(i,j)\in A \suchthat \delta^-(i)=\emptyset,~\delta^+(j) \neq \emptyset \right\},~A_2 \coloneqq \left\{(i,j)\in A \suchthat \delta^-(i) \neq \emptyset,~ \delta^+(j) = \emptyset\right\}, \\
&A_3 \coloneqq \left\{(i,j)\in A \suchthat \delta^-(i) \neq \emptyset,~ \delta^+(j) \neq \emptyset\right\},~A_4\coloneqq \left\{(i,j)\in A \suchthat \delta^-(i)= \emptyset,~ \delta^+(j) = \emptyset\right\},
\end{align*}
and let $\vc{\theta}^\xi = (\vc{\beta}^\xi,\vc{\gamma}^\xi,\vc{\delta}^\xi,\vc{\phi}^\xi,\vc{\lambda}^\xi,\vc{\mu}^\xi)$ be the vector of dual variables associated with constraints of \eqref{obj}-\eqref{x_bound} for a scenario $\xi \in \Xi$.
Further, define the bi-function

\begin{align*}
 f(\vc{y};\vc{\theta}^\xi) = &\sum_{q\in V}\sum_{j\in \delta^+(q)}\left(-\ell_{qj}\beta^\xi_{qj}+u_{qj}\gamma^\xi_{qj}\right) 
+ \sum_{q\in \bar V}\sum_{(i,j)\in \delta^-(q)\times \delta^+(q)}\left(u_{iq}(1-y^q_{ij})\lambda^\xi_{iqj} + u_{qj}(1-y^q_{ij})\mu^\xi_{iqj}\right) \nonumber \\
&+ \sum_{q\in \bar V}\sum_{i\in \delta^-(q)}\left(u_{iq}\sum_{j\in \delta^+(q)}y^q_{ij}\sigma^\xi_{iq}\right) + \sum_{q\in \bar V}\sum_{j\in \delta^+(q)}\left(u_{qj}\sum_{i\in \delta^-(q)}y^q_{ij}\phi^\xi_{qj}\right).
\end{align*}

For a given $\bar{\vc{y}}$ and each scenario $\xi \in \Xi$, the dual of the subproblem~\eqref{obj}-\eqref{x_bound} can be written as follows where the symbol $\star$ on a node means that it belongs to $\bar V$.

\begin{subequations}
\small{
\begin{align}
\min \quad &f(\bar{\vc{y}};\vc{\theta}^\xi) \\
\text{s.t.} \quad &\alpha^\xi_{\accentset{\star}{q}} - \beta^\xi_{i\accentset{\star}{q}} + \gamma^\xi_{i\accentset{\star}{q}} + \sum_{j:j\in\delta^+(\accentset{\star}{q})}\lambda^{\xi}_{i\accentset{\star}{q} j} - \sum_{j:j\in \delta^+(\accentset{\star}{q})}\mu^{\xi}_{i\accentset{\star}{q} j} + \sigma^\xi_{i\accentset{\star}{q}} \ge r_{i\accentset{\star}{q}} &\forall (i,\accentset{\star}{q})\in A_1 \\
&\alpha^\xi_q - \beta^\xi_{iq} + \gamma^\xi_{iq} \ge r_{iq} &\forall (i,q)\in A_1 \\
&-\alpha^\xi_{\accentset{\star}{q}} - \beta^\xi_{\accentset{\star}{q}j} + \gamma^\xi_{\accentset{\star}{q}j} - \sum_{i: i\in \delta^-(\accentset{\star}{q})}\lambda^{\xi}_{i\accentset{\star}{q} j} + \sum_{i: i\in \delta^-(\accentset{\star}{q})}\mu^{\xi}_{i \accentset{\star}{q} j} + \phi^\xi_{\accentset{\star}{q}j} \ge r_{\accentset{\star}{q}j} &\forall (\accentset{\star}{q},j) \in A_2 \\
&-\alpha^\xi_q - \beta^\xi_{qj} + \gamma^\xi_{qj} \ge r_{qj} &\forall (q,j) \in A_2 \\
&-\alpha^\xi_{\accentset{\star}{q}} + \alpha^\xi_{\accentset{\star}{j}} - \beta^\xi_{\accentset{\star}{q}\accentset{\star}{j}} + \gamma^\xi_{\accentset{\star}{q}\accentset{\star}{j}} + \sum_{i\in \delta^-(\accentset{\star}{q})}\left(\mu^{\xi}_{i\accentset{\star}{q} \accentset{\star}{j}}-\lambda^{\xi}_{i\accentset{\star}{q} \accentset{\star}{j}}\right) + \sum_{i\in \delta^+(\accentset{\star}{j})}\left(\lambda^{\xi}_{\accentset{\star}{q} \accentset{\star}{j}i}-\mu^{\xi}_{\accentset{\star}{q} \accentset{\star}{j}i}\right) + \sigma^\xi_{\accentset{\star}{q}\accentset{\star}{j}} + \phi^\xi_{\accentset{\star}{q}\accentset{\star}{j}} \ge r_{\accentset{\star}{q}\accentset{\star}{j}} &\forall (\accentset{\star}{q},\accentset{\star}{j})\in A_3 \\
&-\alpha^\xi_{\accentset{\star}{q}} + \alpha^\xi_j - \beta^\xi_{\accentset{\star}{q}j} + \gamma^\xi_{\accentset{\star}{q}j} + \sum_{i\in \delta^-(\accentset{\star}{q})}\left(\mu^{\xi}_{i\accentset{\star}{q} j}-\lambda^{\xi}_{i\accentset{\star}{q} j}\right) + \phi^\xi_{\accentset{\star}{q}j} \ge r_{\accentset{\star}{q}j} &\forall (\accentset{\star}{q},j)\in A_3 \\
&-\alpha^\xi_q + \alpha^\xi_{\accentset{\star}{j}} - \beta^\xi_{q\accentset{\star}{j}} + \gamma^\xi_{q\accentset{\star}{j}} + \sum_{i\in \delta^+(\accentset{\star}{j})}\left(\lambda^{\xi}_{q\accentset{\star}{j} i}-\mu^{\xi}_{q\accentset{\star}{j} i}\right) + \sigma^\xi_{q\accentset{\star}{j}} \ge r_{q\accentset{\star}{j}} &\forall (q,\accentset{\star}{j})\in A_3 \\
&-\alpha^\xi_q + \alpha^\xi_j - \beta^\xi_{qj} + \gamma^\xi_{qj} \ge r_{qj} &\forall (q,j)\in A_3 \\
&-\beta^\xi_{iq} + \gamma^\xi_{iq} \ge r_{iq} &\forall (i,q)\in A_4 \\
&\alpha^\xi_q \in \mathbb{R} &\forall q\in V' \\  
&\beta^\xi_{ij},~\gamma^\xi_{ij},~\sigma^\xi_{ij},~\phi^\xi_{ij},~\lambda^\xi_{iqj},~\mu^\xi_{iqj} \ge 0 &\forall i,q,j\in V. 
\end{align}
}
\end{subequations}

If the above problem has an optimal solution $\hat{\vc \theta}^\xi$ for all $\xi \in \Xi$, the output of the subproblems will be an optimality cut of the form $\sum_{\xi \in \Xi}\prob^\xi f(\vc{y};\hat{\vc{\theta}}^\xi) \geq z$.
If the above problem is unbounded along a ray $\hat{\vc{\theta}}^\xi$ for a $\xi \in \Xi$, the output of the subproblem will be a feasibility cut of the form $f(\vc{y};\hat{\vc{\theta}}^\xi) \geq 0$.
Note that replacing variables $\vc{y}$ in the above constraints with $\vc{w}$ through the mapping~\eqref{eq:y-w} results in separable nonlinear constraints.
Nevertheless, since these constraints will be used to refine the master DD, their incorporation is simple due to structural flexibility of DDs in modeling such constraints; we refer the reader to \cite{davarnia2020outer} for a detailed account for modeling INLPs with DDs.

\section{Computational Experiments} \label{sec: Computational Experiments}
In this section, we solve SGUFP as a core model for the unit train scheduling problem with demand stochasticity using three different approaches: (i) the standard MIP formulation that is a deterministic equivalent of the two-stage model and contains all variables and constraints of the master problem and $|\Xi|$ subproblems; (ii) the Benders reformulation presented in Section~\ref{subsec: MIP Formulation} composed of the master problem~\eqref{obj1}-\eqref{binary_y} and $|\Xi|$ subproblems~\eqref{obj}-\eqref{x_bound}; and (iii) the DD-BD algorithm proposed in the present paper. 
In the Benders approach, we solve separate subproblems using a fixed vector $\bar{\vc{y}}$ obtained from the master problem. The feasibility cuts generated by subproblems are added directly to the constraint set of the master problem, and the optimality cuts are added as an aggregated cut over all scenarios.
We note here that when there is a feasibility cut for any scenario, we add it directly to separate the solution of the current iteration and move on to the next
iteration.
To obtain a valid inequality that provides a bound for the single $z$ variable, we need to aggregate valid inequalities over all scenario subproblems as $z$ is composed of the objective value of all these subproblems. Therefore, we can only produce an optimality cut for the $z$ variable when we have optimality cuts for all of the subproblems.
For the DD-BD approach, we use the following algorithmic choices to build restricted and relaxed DDs.
For the restricted DDs, we choose a subset of the $r$-$t$ paths with largest lengths, which are more likely to contain an optimal solution.
For the relaxed DDs, we merge nodes that have the largest number of common members in their state values.
We refer the reader to \cite{bergman2016decision} for other heuristic approaches that could be used for this purpose.

\subsection{Test Instances} \label{subsec: Testbeds}
In our experiments, we consider the structure of the SGUFP network given in Section~\ref{subsec: MIP Formulation}.
To ensure that the problem is always feasible, we create an artificial node $s_0$ to compensate for any shortage of the supply, and add an arc from the artificial supply $s_0$ to each demand node.

We create test instances based on the specification given in \cite{davarnia2019network}, which is inspired by realistic models.
In particular, we consider a base rail network $G'=(V',A')$ where $10\%$ and $30\%$ of the nodes are supply and demand nodes, respectively.
We assume that $50\%$ of the nodes must satisfy the NSNM requirement. 
We then create a network $G=(V,A)$ by augmenting supply/demand and artificial nodes as described above with the following settings. 
The integer supply value at supply nodes is randomly selected from the interval $[100,600]$. 
The capacity of arcs connecting $s_0$ to demand nodes are considered to be unbounded, and the integer capacity value of other arcs is randomly selected from the interval $[100,300]$. 
For each demand scenario $\xi \in \Xi$, the integer demand value at demand nodes is randomly chosen from the interval $[100,200]$. 
The reward of the arcs connecting $s_0$ to the demand nodes are generated from the interval $[-10,-5]$ to represent the cost of lost demands. 
The reward of the arcs connecting the source to the supply nodes is randomly selected from the interval $[5,10]$, and the reward of the arcs connecting the demand nodes to the sink is fixed to zero since the flow of these arcs is also fixed. 
The reward of all other arcs is created randomly from the interval $[-2,2]$ where the negative values indicate the cost of sending flows through congested arcs.  
We consider four categories of rail networks with $|V'|\in\{40,60,80,100\}$. 
For each category, we create five scenario classes for the number of demand scenarios $|\Xi|\in\{50,100,150,200,250\}$.
For each network category and scenario class, we create five random instances based on the above settings. Test instances are publicly available~\citep{salemi2022}.
\subsection{Numerical Results}
\begin{figure}[hbt!]
\centering
\includegraphics[scale=0.80]{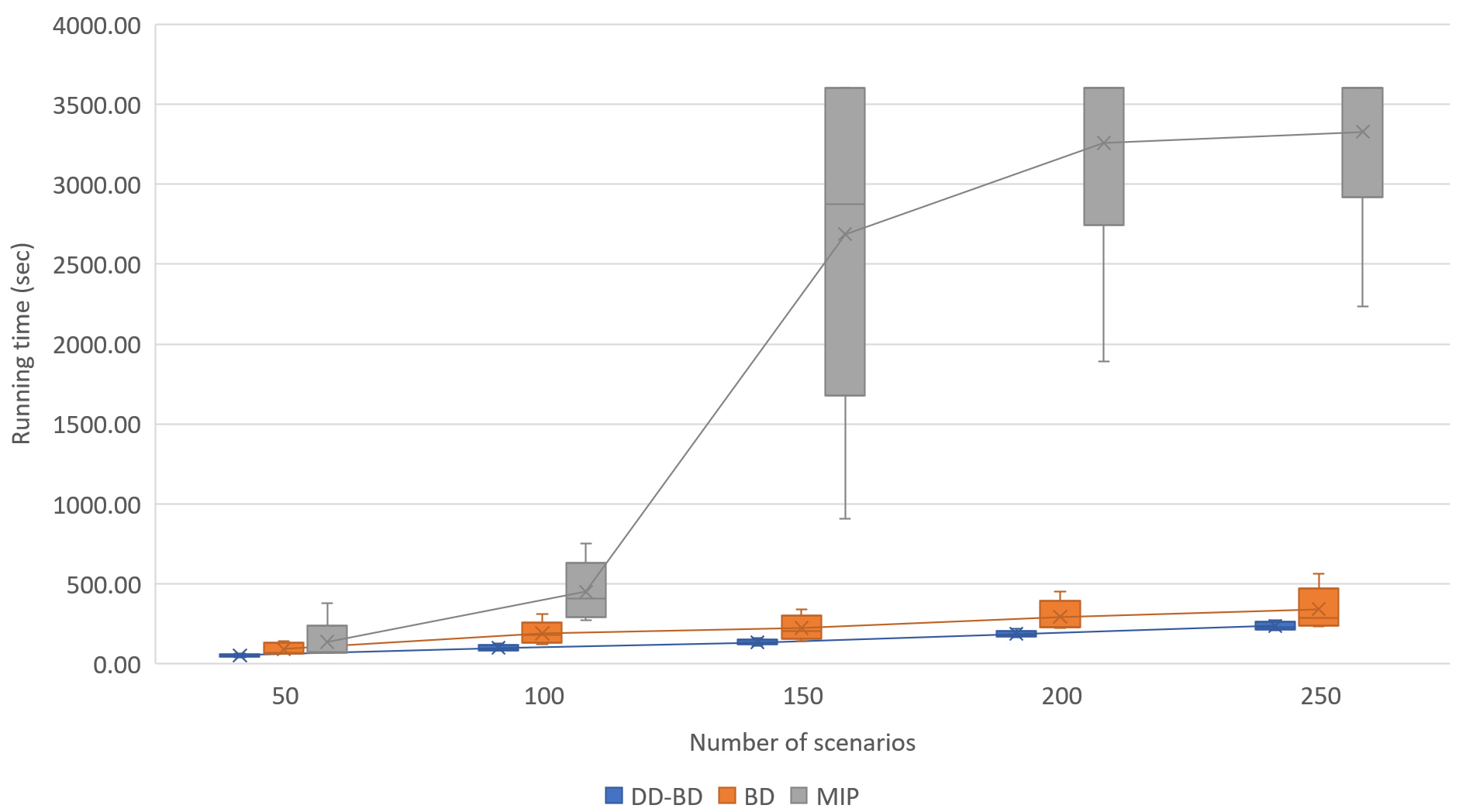}
\caption{Comparison of DD-BD, BD, and MIP models when $|V'|=40$ under five scenarios}
\label{fig: n40}
\end{figure}

\begin{figure}[hbt!]
\centering
\includegraphics[scale=0.80]{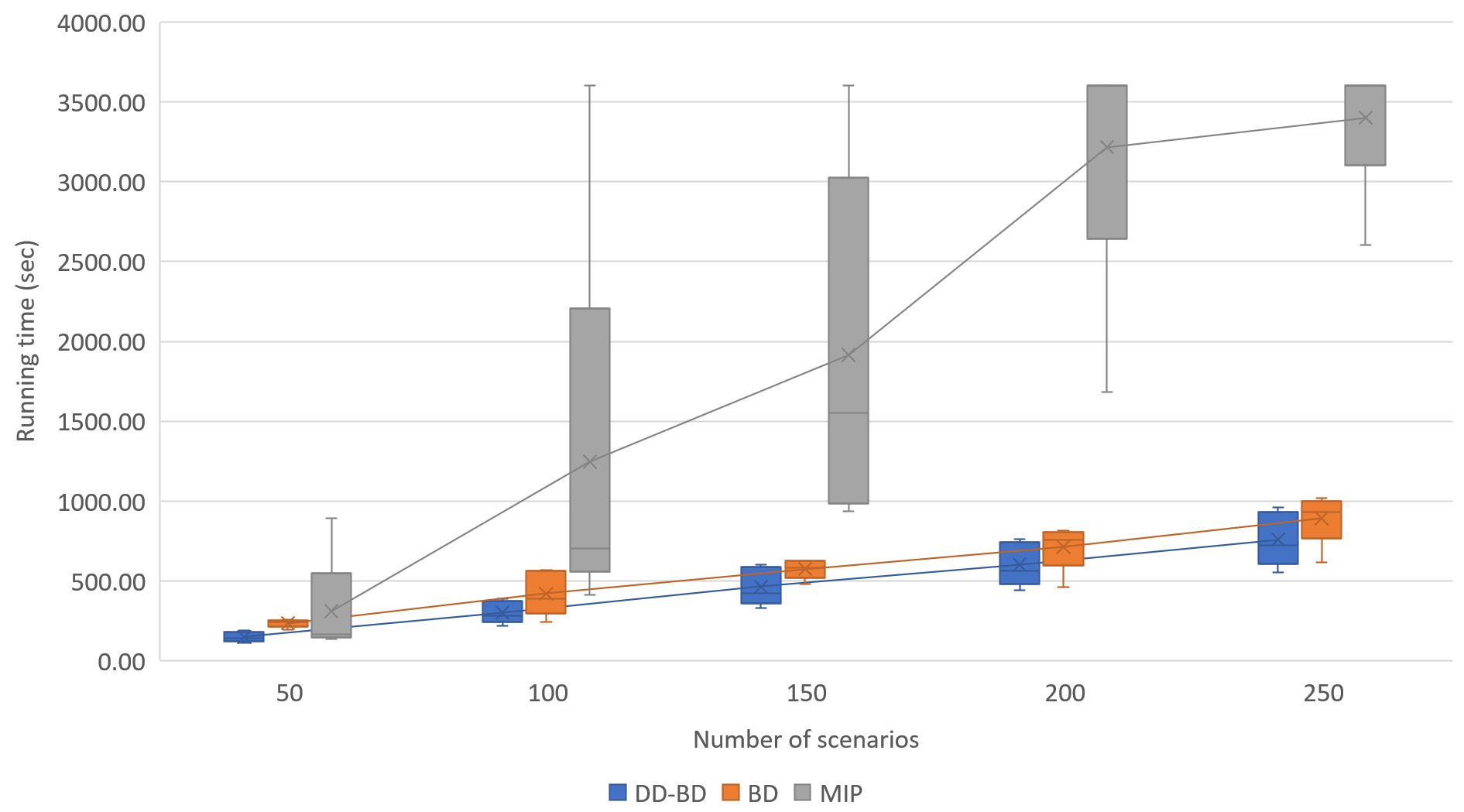}
\caption{Comparison of DD-BD, BD, and MIP models when $|V'|=60$ under five scenarios}
\label{fig: n60}
\end{figure}

\begin{figure}[hbt!]
\centering
\includegraphics[scale=0.80]{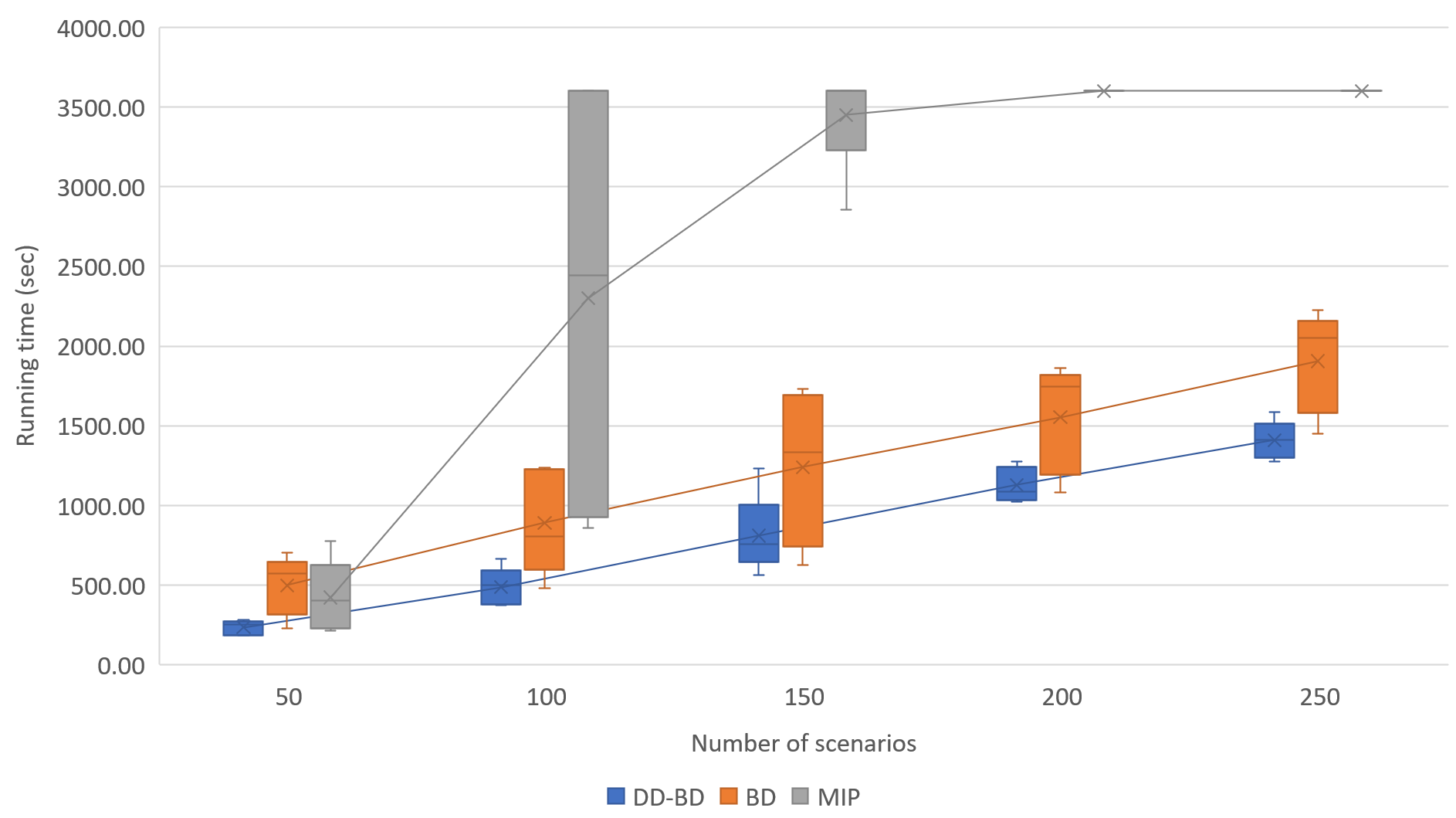}
\caption{Comparison of DD-BD, BD, and MIP models when $|V'|=80$ under five scenarios}
\label{fig: n80}
\end{figure}

\begin{figure}[hbt!]
\centering
\includegraphics[scale=0.80]{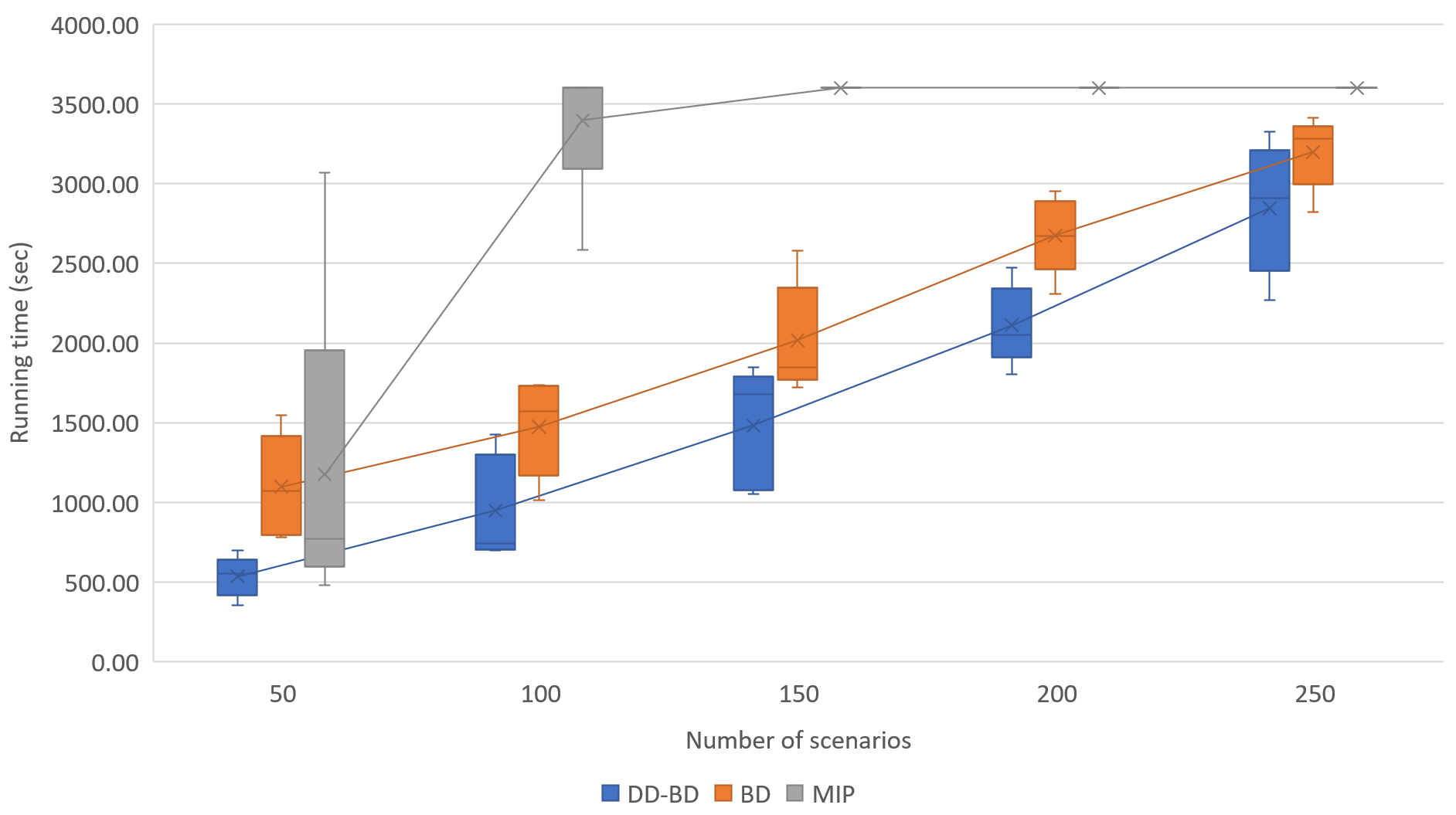}
\caption{Comparison of DD-BD, BD, and MIP models when $|V'|=100$ under five scenarios}
\label{fig: n100}
\end{figure}

In this section, we present the numerical results that compare the performance of the DD-BD formulation for the SGUFP instances with that of the MIP formulation, denoted by ``MIP", and the standard Benders reformulation, denoted by ``BD". All experiments are conducted on a machine running Windows 10, x64 operating system with Intel\textsuperscript{\textregistered} Core i7 processor (2.60 GHz) and 32 GB RAM.  
The Gurobi optimization solver (version 9.1.1) is used to solve instances for the MIP and BD models. When solving problems with Gurobi, we turn off presolve and cuts for all methods to have a fair comparison. Tables~\ref{tab:1}-\ref{tab:4} report the running times of each of these formulations for $|V'|\in\{40,60,80,100\}$ and $|\Xi|\in\{50,100,150,200,250\}$ where the time limit is set to 3600 seconds. 
The symbol $``>3600"$ indicates that the problem was not solved within the time limit.
As evident in these tables, the DD-BD formulation outperforms the other alternatives.
In particular, the gap between the solution time of the DD-BD and the MIP and BD approaches widens as the problem size increases.
For example, as reported in Table~\ref{tab:1}, while the DD-BD approach solves all 25 instances in under 275 seconds, the MIP approach fails to solve 10 of them within 3600 seconds, $80\%$ of which involve 200 or 250 scenarios. This shows a clear superiority of the DD-BD over the MIP method.
Further, for most of the instances, the DD-BD approach outperforms the standard BD approach, rendering it as the superior solution method among all three. 
Figures~\ref{fig: n40}-\ref{fig: n100} compare the performance of DD-BD, BD, and MIP formulations through box and whisker plots for each network size and under each scenario class. 
In these figures, for uniformity of illustration, we used 3600 seconds for the running time of instances that fail to solve the problem within that time limit.  
As the figures show, the minimum, median, and maximum of running times of the DD-BD method are remarkably smaller than those of the both BD and MIP methods in all cases. 
These results show the potential of the DD-BD framework in solving network problems with challenging combinatorial structures.  
In Appendix~\ref{app:limit}, we present additional numerical results for the DD-BD approach to assess its ability to solve larger problem sizes.

\begin{table}[!htbp]
\centering
\caption{Running times (in seconds) of MIP, BD, and DD-BD for $|V'|=40$.}
\label{tab:1}
\scalebox{0.9}{
\begin{tabular}{c|l|rrrrr}
\multirow{2}{*}{Instance \#} & \multirow{2}{*}{Model} & \multicolumn{5}{c}{Number of scenarios} \\
 &  & 50 & 100 & 150 & 200 & 250 \\\hline
\multirow{3}{*}{1} 
 & MIP & 75.74 & 512.62 & 2877.19 & $>3600$ & $>3600$ \\
 & BD & 141.83 & 313.84 & 339.81 & 451.93 & 565.82 \\
 & DD-BD & 56.94 & 129.87 & 163.43 & 219.02 & 274.36 \\\hline 
 \multirow{3}{*}{2} 
 & MIP & 67.59 & 275.07 & 906.10 & 1892.21 & 2235.53 \\
 & BD & 63.44 & 121.25 & 141.04 & 230.81 & 235.87 \\
 & DD-BD & 42.60 & 82.65 & 128.16 & 164.52 & 208.94 \\\hline 
 \multirow{3}{*}{3} 
 & MIP & 94.86 & 753.23 & 2453.05 & $>3600$ & $>3600$ \\
 & BD & 71.14 & 139.20 & 172.86 & 224.33 & 244.91 \\
 & DD-BD & 53.32 & 93.58 & 113.93 & 178.65 & 217.33 \\\hline 
 \multirow{3}{*}{4} 
 & MIP & 71.46 & 309.62 & $>3600$ & $>3600$ & $>3600$ \\
 & BD & 63.55 & 182.01 & 267.94 & 334.74 & 380.22 \\
 & DD-BD & 46.61 & 87.81 & 130.19 & 183.23 & 253.72 \\\hline 
 \multirow{3}{*}{5} 
 & MIP & 380.33 & 406.73 & $>3600$ & $>3600$ & $>3600$ \\
 & BD & 123.69 & 198.73 & 205.16 & 231.56 & 287.24 \\
 & DD-BD & 67.04 & 104.78 & 138.46 & 195.69 & 231.74 \\
\end{tabular}}
\end{table} 

\begin{table}[!htbp]
\centering
\caption{Running times (in seconds) of MIP, BD, and DD-BD for $|V'|=60$.}
\label{tab:2}
\scalebox{0.9}{
\begin{tabular}{c|l|rrrrr}
\multirow{2}{*}{Instance \#} & \multirow{2}{*}{Model} & \multicolumn{5}{c}{Number of scenarios} \\
 &  & 50 & 100 & 150 & 200 & 250 \\\hline
\multirow{3}{*}{1} 
 & MIP & 893.73 & $>3600$ & $>3600$ & $>3600$ & $>3600$ \\
 & BD & 241.85 & 556.18 & 582.80 & 758.54 & 933.05 \\
 & DD-BD & 176.16 & 357.06 & 603.81 & 719.27 & 901.02 \\\hline 
 \multirow{3}{*}{2} 
 & MIP & 206.87 & 811.64 & 1554.10 & $>3600$ & $>3600$ \\
 & BD & 259.63 & 351.39 & 624.08 & 816.44 & 1017.95 \\
 & DD-BD & 189.07 & 388.85 & 572.52 & 764.76 & 961.35 \\\hline 
 \multirow{3}{*}{3} 
 & MIP & 139.70 & 702.96 & 1035.79 & $>3600$ & $>3600$ \\
 & BD & 246.48 & 569.37 & 628.84 & 795.56 & 978.15 \\
 & DD-BD & 142.81 & 284.65 & 422.52 & 565.23 & 725.86 \\\hline 
 \multirow{3}{*}{4} 
 & MIP & 153.16 & 415.46 & 938.03 & 1681.21 & 2604.25 \\
 & BD & 238.33 & 388.19 & 563.15 & 732.59 & 919.08 \\
 & DD-BD & 131.29 & 262.36 & 393.18 & 521.12 & 654.71 \\\hline 
 \multirow{3}{*}{5} 
 & MIP & 165.57 & 706.16 & 2447.15 & $>3600$ & $>3600$ \\
 & BD & 194.12 & 244.61 & 479.32 & 463.63 & 617.09 \\
 & DD-BD & 112.09 & 221.30 & 332.25 & 443.96 & 556.33 \\
\end{tabular}}
\end{table}

\begin{table}[!htbp]
\centering
\caption{Running times (in seconds) of MIP, BD, and DD-BD for $|V'|=80$.}
\label{tab:3}
\scalebox{0.9}{
\begin{tabular}{c|l|rrrrr}
\multirow{2}{*}{Instance \#} & \multirow{2}{*}{Model} & \multicolumn{5}{c}{Number of scenarios} \\
 &  & 50 & 100 & 150 & 200 & 250 \\\hline
\multirow{3}{*}{1} 
 & MIP & 215.82 & 860.21 & $>3600$ & $>3600$ & $>3600$ \\
 & BD & 588.51 & 806.61 & 1731.50 & 1860.12 & 2051.52 \\
 & DD-BD & 256.12 & 500.52 & 757.68 & 1025.88 & 1278.13 \\\hline 
 \multirow{3}{*}{2} 
 & MIP & 479.76 & $>3600$ & $>3600$ & $>3600$ & $>3600$ \\
 & BD & 398.29 & 713.01 & 861.65 & 1080.79 & 1709.04 \\
 & DD-BD & 184.34 & 379.04 & 724.66 & 1088.21 & 1587.90 \\\hline 
 \multirow{3}{*}{3} 
 & MIP & 238.79 & 996.22 & $>3600$ & $>3600$ & $>3600$ \\
 & BD & 702.18 & 1236.58 & 1650.42 & 1773.63 & 2227.89 \\
 & DD-BD & 285.13 & 518.46 & 778.97 & 1046.39 & 1326.22 \\\hline 
 \multirow{3}{*}{4} 
 & MIP & 404.26 & 2441.64 & 2855.29 & $>3600$ & $>3600$ \\
 & BD & 572.83 & 1219.37 & 1334.21 & 1745.91 & 2089.80 \\
 & DD-BD & 263.78 & 665.30 & 1230.81 & 1277.93 & 1444.02  \\\hline 
 \multirow{3}{*}{5} 
 & MIP & 778.50 & $>3600$ & $>3600$ & $>3600$ & $>3600$ \\
 & BD & 231.11 & 481.31 & 625.91 & 1310.24 & 1452.27 \\
 & DD-BD & 187.34 & 376.96 & 564.34 & 1205.54 & 1412.94 \\ 
\end{tabular}}
\end{table}

\begin{table}[!htbp]
\centering
\caption{Running times (in seconds) of MIP, BD, and DD-BD for $|V'|=100$.}
\label{tab:4}
\scalebox{0.9}{
\begin{tabular}{c|l|rrrrr}
\multirow{2}{*}{Instance \#} & \multirow{2}{*}{Model} & \multicolumn{5}{c}{Number of scenarios} \\
 &  & 50 & 100 & 150 & 200 & 250 \\\hline
\multirow{3}{*}{1} 
 & MIP & 774.18 & $>3600$ & $>3600$ & $>3600$ & $>3600$ \\
 & BD & 1282.59 & 1728.71 & 1848.49 & 2307.74 & 3309.93 \\
 & DD-BD & 698.36 & 1427.38 & 1731.95 & 2014.96 & 3323.54 \\\hline 
 \multirow{3}{*}{2} 
 & MIP & 480.97 & $>3600$ & $>3600$ & $>3600$ & $>3600$ \\
 & BD & 781.47 & 1573.23 & 1820.79 & 2672.18 & 2819.61 \\
 & DD-BD & 586.89 & 1171.96 & 1848.49 & 2471.49 & 2635.22 \\\hline 
 \multirow{3}{*}{3} 
 & MIP & 3071.37 & $>3600$ & $>3600$ & $>3600$ & $>3600$ \\
 & BD & 1072.14 & 1322.96 & 2112.50 & 2951.55 & 3412.99 \\
 & DD-BD & 485.31 & 703.70 & 1055.36 & 1803.66 & 2269.97 \\\hline 
 \multirow{3}{*}{4} 
 & MIP & 838.79 & 2585.38 & $>3600$ & $>3600$ & $>3600$ \\
 & BD & 1548.93 & 1738.92 & 2580.53 & 2616.19 & 3169.28 \\
 & DD-BD & 554.89 & 743.64 & 1098.82 & 2052.73 & 3094.23 \\\hline 
 \multirow{3}{*}{5} 
 & MIP & 714.39 & $>3600$ & $>3600$ & $>3600$ & $>3600$ \\
 & BD & 808.48 & 1013.68 & 1722.01 & 2824.14 & 3282.10 \\
 & DD-BD & 353.48 & 700.57 & 1680.60 & 2213.81 & 2907.78 \\
\end{tabular}}
\end{table}

We conclude this section by noting that, while the focus of this paper has been on the unit train problem with the no-split no-merge requirements, the proposed DD-BD framework can be applied to model network problems that contain additional side constraints on the flow variables, as those constraints can be handled in the subproblems while the DD structure in the master problem remains intact.
Examples of such side constraints include the \emph{usage-fee} limitation~\citep{holzhauser2017network} and the \emph{flow ratio} requirement~\citep{holzhauser2017maximum}. 
Applying the DD-BD method to such network models and assessing its effectiveness compared to alternative approaches could be an interesting direction for future research.

\section{Conclusion} \label{sec: Conclusion}
In this paper, we introduce a DD-based framework to solve the SGUFP. This framework uses Benders decomposition to decompose the SGUFP into a master problem composed of the combinatorial NSNM constraints, and a subproblem that solves a continuous network flow model. The master problem is modeled by a DD, which is successively refined with respect to the cuts generated through subproblems. To assess the performance of the proposed method, we apply it to a variant of unit train scheduling problem formulated as a SGUFP, and compare it with the standard MIP and Benders reformulation of the problem.

\ACKNOWLEDGMENT{This project is sponsored in part by the Iowa Energy Center, Iowa Economic Development Authority
	and its utility partners. We thank the anonymous referees and the Associate Editor for their helpful comments that contributed to improving the paper.
}


\clearpage
\bibliographystyle{informs2014trsc} 
\bibliography{bib} 


%
%

\clearpage

\begin{APPENDICES}
\section{Comparison of Master Problem Formulations} \label{app:comparison}

In this section, we describe the differences between DDs in the space of $\vc{w}$ variables and those in the space of original $\vc{y}$ in the master problem formulation \eqref{eq:master} in Section~\ref{subsec:master}.
First, we illustrate the size difference between these DDs in Example~\ref{ex:DD-w}.

\begin{example} \label{ex:DD-w}
Consider a directed graph $G=(V,A)$ with node set $V = \{1,2,q,3,4\}$ and arc set $A=\{(1,q),(2,q),(q,3),(q,4)\}$ where the central node $q$ is subject to NSNM constraints. Let $\ind^-(1,q) = \ind^+(q,3) = 1$ and $\ind^-(2,q) = \ind^+(q,4) = 2$. Then, the exact DDs showed in Figures~\ref{subfig: DD_1} and~\ref{subfig: DD_2} with three and five arc layers represent the feasible region of master problem~\eqref{eq:master} and~\eqref{obj1}-\eqref{binary_y}, respectively, where $-M$ and $M$ are valid bounds for variable $z$.

\begin{figure}[!htbp]
\captionsetup[subfigure]{justification=centering}
\centering
\begin{subfigure}[b]{0.45\linewidth} 
\centering
\includegraphics[scale=0.25]{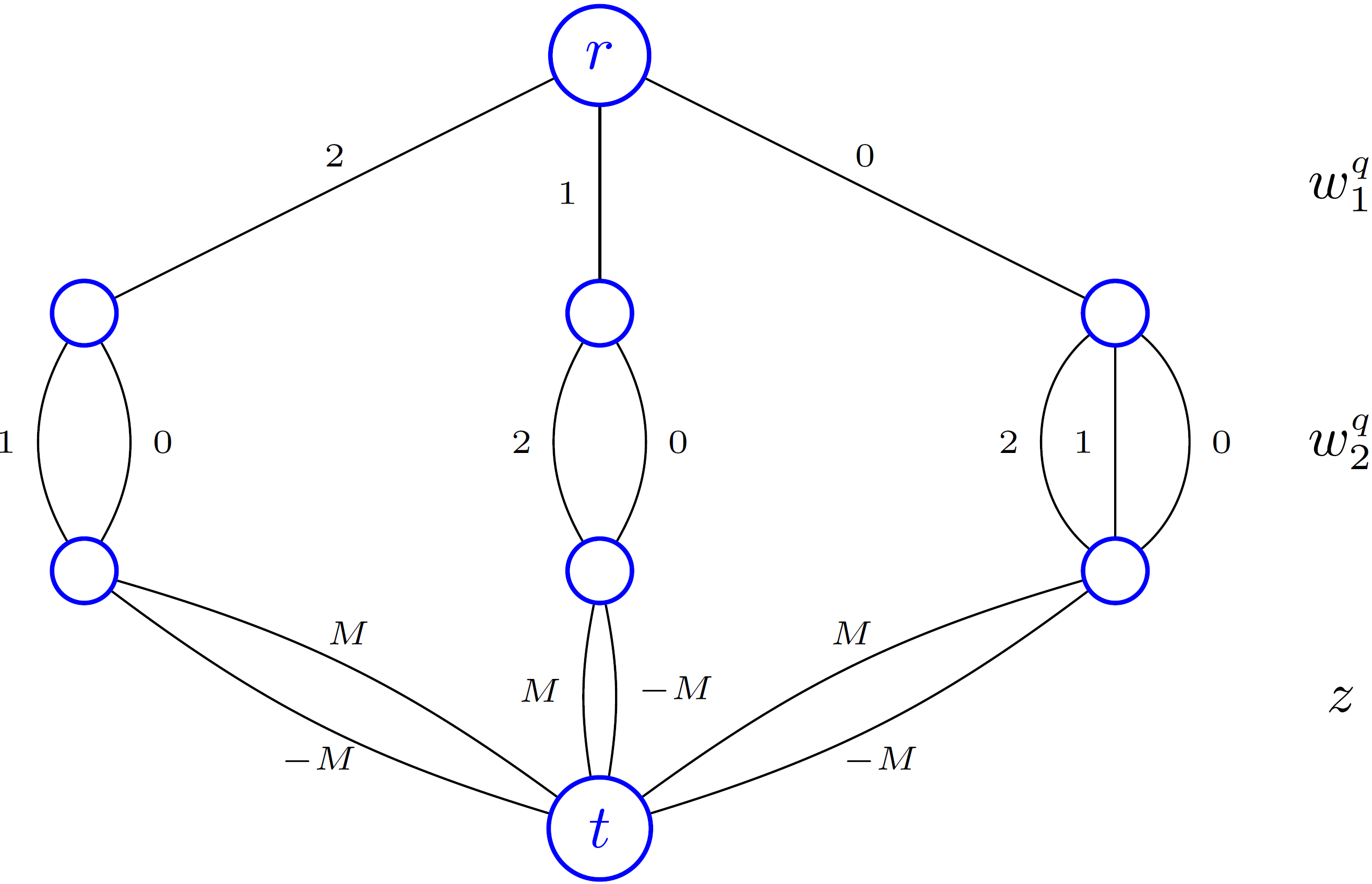}
\caption{A DD in the space of $\vc w$ variables. Numbers next to arcs represent labels.}
\label{subfig: DD_1}
\end{subfigure}\hspace{1cm}
\begin{subfigure}[b]{0.45\linewidth} 
\centering
\includegraphics[scale=0.25]{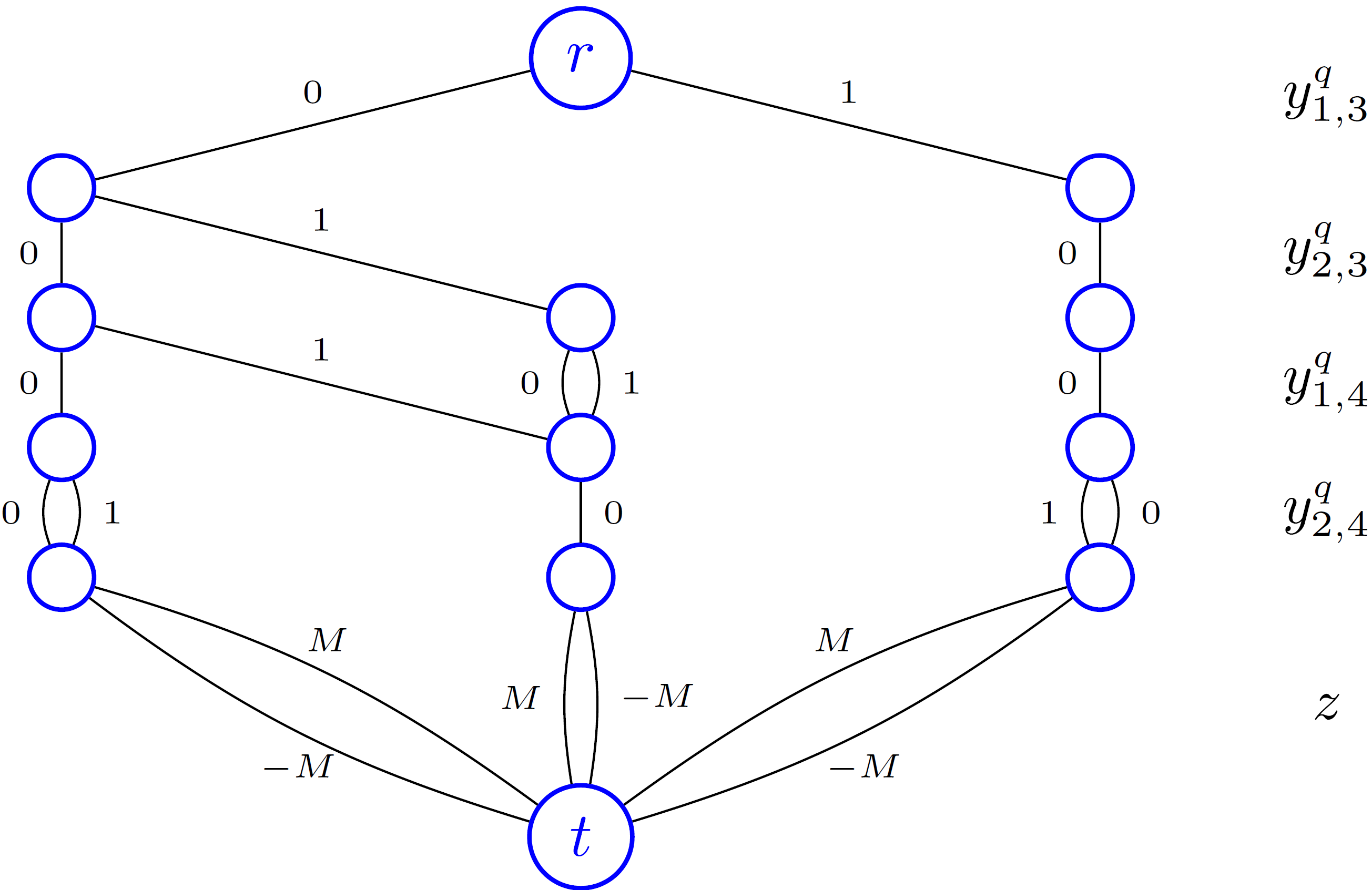}
\caption{A DD in the space of $\vc y$ variables. Numbers next to arcs represent labels.}
\label{subfig: DD_2}
\end{subfigure}
\caption{Comparison of the number of arc layers for DDs in the space of $\vc w$ and $\vc y$ variables}
\end{figure}

\end{example}

As evident from the above example, the main advantage of using a DD in the space of $\vc{w}$ is the reduction in the number of arc layers, which is the main determinant of the DDs computational efficiency. 
In particular, even though such a DD has a larger number of nodes at the layers, a relaxed DD can be constructed to limit the width, and hence provide an efficient relaxed DD in a smaller dimension, whereas the relaxations of the DD constructed in the space of $\vc{y}$ variables would still be higher-dimensional.

To assess the computational efficiency of the solution approach in relation to the DD space, we compare the performance of the DD-BD method under two different settings: (i) where DDs are built in the space of $\vc w$ variables, denoted by DD-BD-$\vc w$, and (ii) where DDs are built in the space of $\vc y$ variables, denoted by DD-BD-$\vc y$. We report the results of these two implementations for $|V'| \in \{40,80\}$ and under five different scenarios in Table~\ref{tab:DDs} and Table~\ref{tab:DDs-80}.

As observed in these tables, the DD-BD-$\vc w$ solves all instances faster than DD-BD-$\vc y$, with orders of magnitude time improvement as the problem size (number of scenarios) increases.
These preliminary computational results show the advantage of designing the DD-BD method for the SGUFP in a transformed space of variables.

\begin{table}[!htbp]
\centering
\caption{Running times (in seconds) of DD-BD-$\vc w$ and DD-BD-$\vc y$ for $|V'|=40$.}
\label{tab:DDs}
\scalebox{0.9}{
\begin{tabular}{c|l|rrrrr}
\multirow{2}{*}{Instance \#} & \multirow{2}{*}{Model} & \multicolumn{5}{c}{Number of scenarios} \\
&  & 50 & 100 & 150 & 200 & 250 \\\hline
\multirow{2}{*}{1} 
 & DD-BD-$\vc w$ & 56.94 & 129.87 & 163.43 & 219.02 & 274.36 \\
 & DD-BD-$\vc y$ & 89.68  & 304.08 & 432.34 & 642.70 & 839.57 \\\hline 
 \multirow{2}{*}{2} 
 & DD-BD-$\vc w$ & 42.60 & 82.65 & 128.16 & 164.52 & 208.94 \\
 & DD-BD-$\vc y$ & 68.23 & 148.76 & 244.53 & 344.86 & 605.04 \\\hline 
 \multirow{2}{*}{3} 
 & DD-BD-$\vc w$ & 53.32 & 93.58 & 113.93 & 178.65 & 217.33 \\
 & DD-BD-$\vc y$ & 83.05 & 157.67 & 310.07 & 541.33 & 658.98 \\\hline 
 \multirow{2}{*}{4} 
 & DD-BD-$\vc w$ & 46.61 & 87.81 & 130.19 & 183.23 & 253.72 \\
 & DD-BD-$\vc y$ & 78.11 & 149.26 & 325.31 & 460.73 & 694.57 \\\hline 
 \multirow{2}{*}{5} 
 & DD-BD-$\vc w$ & 67.04 & 104.78 & 138.46 & 195.69 & 231.74 \\
 & DD-BD-$\vc y$ & 109.61 & 223.78 & 351.80 & 532.12 & 669.78 
\end{tabular}}
\end{table} 

\begin{table}[!htbp]
\centering
\caption{Running times (in seconds) of DD-BD-$\vc w$ and DD-BD-$\vc y$ for $|V'|=80$.}
\label{tab:DDs-80}
\scalebox{0.9}{
\begin{tabular}{c|l|rrrrr}
\multirow{2}{*}{Instance \#} & \multirow{2}{*}{Model} & \multicolumn{5}{c}{Number of scenarios} \\
&  & 50 & 100 & 150 & 200 & 250 \\\hline
\multirow{2}{*}{1} 
 & DD-BD-$\vc w$ & 256.12 & 500.52 & 757.68 & 1025.88 & 1278.13 \\
 & DD-BD-$\vc y$ & 483.42 & 977.03 & 1642.27 & 3175.72 & 4230.29 \\\hline 
 \multirow{2}{*}{2} 
 & DD-BD-$\vc w$ & 184.34 & 379.04 & 724.66 & 1088.21 & 1587.90 \\
 & DD-BD-$\vc y$ & 340.13 & 864.21 & 1856.96 & 3010.55 & 4843.67   \\\hline 
 \multirow{2}{*}{3} 
 & DD-BD-$\vc w$ & 285.13 & 518.46 & 778.97 & 1046.39 & 1326.22 \\
 & DD-BD-$\vc y$ & 568.32 & 1176.44 & 2401.98 & 3326.76 & 4283.58  \\\hline 
 \multirow{2}{*}{4} 
 & DD-BD-$\vc w$ & 263.78 & 665.30 & 1230.81 & 1277.93 & 1444.02 \\
 & DD-BD-$\vc y$ & 501.04 & 1430.77 & 2868.92 & 3356.39 & 4356.48  \\\hline 
 \multirow{2}{*}{5} 
 & DD-BD-$\vc w$ & 187.34 & 376.96 & 564.34 & 1205.54 & 1412.94 \\
 & DD-BD-$\vc y$ & 354.37 & 781.18 & 1279.73 & 3001.72 & 3834.08  
\end{tabular}}
\end{table}

\section{Additional Computational Experiments} \label{app:limit}

In this section, we present additional numerical results to assess the limits of the DD-BD method for larger problem instances.
These results are given in Tables~\ref{tab:A1} and~\ref{tab:A2}, where the columns are defined similarly to those of Tables 1-4.
For these instances, the time limit is set to 3600 seconds, and the symbol “$> 3600$” indicates that the problem is not
solved within this time limit.

\begin{table}[!htbp]
\centering
\caption{Running times (in seconds) of DD-BD for $|V'|=120$.}
\label{tab:A1}
\scalebox{1}{
\begin{tabular}{c|l|rrrrr}
\multirow{2}{*}{Instance \#} & \multirow{2}{*}{Model} & \multicolumn{5}{c}{Number of scenarios} \\
 &  & 50 & 100 & 150 & 200 & 250 \\\hline
1 & DD-BD & 1494.49 & 2824.58 & $>3600$ & $>3600$ & $>3600$ \\
2 & DD-BD & 975.47 & 1892.41 & 3198.18 & $>3600$ & $>3600$ \\
3 & DD-BD & 1150.30 & 2263.09 & 3454.47 & $>3600$ & $>3600$ \\
4 & DD-BD & 1261.59 & 2403.79 & $>3600$ & $>3600$ & $>3600$ \\
5 & DD-BD & 906.34 & 1863.15 & 3050.68 & $>3600$ & $>3600$
\end{tabular}}
\end{table}

\begin{table}[!htbp]
\centering
\caption{Running times (in seconds) of DD-BD for $|V'|=150$.}
\label{tab:A2}
\scalebox{1}{
\begin{tabular}{c|l|rrrrr}
\multirow{2}{*}{Instance \#} & \multirow{2}{*}{Model} & \multicolumn{5}{c}{Number of scenarios} \\
 &  & 50 & 100 & 150 & 200 & 250 \\\hline
1 & DD-BD & 2496.16 & $>3600$ & $>3600$ & $>3600$ & $>3600$ \\
2 & DD-BD & 2944.20 & $>3600$ & $>3600$ & $>3600$ & $>3600$ \\
3 & DD-BD & 2321.62 & $>3600$ & $>3600$ & $>3600$ & $>3600$ \\
4 & DD-BD & 2590.34 & $>3600$ & $>3600$ & $>3600$ & $>3600$ \\
5 & DD-BD & 2298.36 & $>3600$ & $>3600$ & $>3600$ & $>3600$
\end{tabular}}
\end{table}

\end{APPENDICES}

\end{document}